\documentclass[10pt]{amsart}
\usepackage{amsmath, latexsym, amssymb, lscape, amsthm}
\usepackage[usenames]{color}
\usepackage[all]{xy}
\input{grcawol.sty}
\newcommand{\gwnot}[2]{
   \grcalca = \grcolumn
   \multiply \grcalca by \factor
   \grcalcc = #1
   \multiply \grcalcc by \hfactor
   \advance \grcalca by \grcalcc
   \grcalcb = \grrow
   \multiply \grcalcb by \factor
   \advance \grcalcb by -\hfactor
   \put(\grcalca,\grcalcb) {\makebox(0,0){$\scriptstyle #2$}} }
\newcommand\gdnot[1]{\gwnot2{#1}}
\def\BC{\mathbb{C}}

\def\qed{{\ \ \ \mbox{$\square$}}}

\def\lra{{\longrightarrow}}

\def\End{{\mbox{\bf{End}}}}

\newcommand\C[1]{{#1\mbox{-\bf{mod}}_{\operatorname{\mathsf fin}}}}
\newcommand\FPdim{\operatorname{FPdim}}

\newtheorem{thm}{Theorem}[section]
\newtheorem{cor}[thm]{Corollary}
\newtheorem{defn}[thm]{Definition}
\newtheorem{prop}[thm]{Proposition}
\newtheorem{remark}[thm]{Remark}

\newtheorem{lem}[thm]{Lemma}

\newcommand\leaveout[1]{}
\newcommand\ptr{\operatorname{\underline{ptr}}}
\newcommand\ptrl{\ptr^\ell}
\newcommand\ptrr{\ptr^r}
\newcommand\FS{\operatorname{\underline{FS}}}
\newcommand\Unit{U}
\newcommand\Counit{C}
\newcommand\gsym{\gcn1113\gcn111{-1}}
\newcommand\triv{{\operatorname{triv}}}
\newcommand\Triv{{\underline{\operatorname{Triv}}}}
\newcommand\redu[1]{{^\vee#1}}
\newcommand\sym{{\operatorname{sym}}}

\newcommand\id{\operatorname{id}}
\renewcommand\o{\otimes}
\newcommand\ro{\otimes}
\newcommand\Tr{\operatorname{Tr}}

\newcommand\piv{j}
\newcommand\dual[1]{#1^\vee}
\newcommand\inv{^{-1}}
\DeclareMathOperator\ev{{\operatorname{ev}}}
\DeclareMathOperator\db{\operatorname{db}}
\newcommand\CC{\mathcal C}
\newcommand\DD{\mathcal D}
\newcommand\FF{\mathcal F}
\newcommand\bidual[1]{#1^{\vee\vee}}
\newcommand\du{^{\vee}}
\newcommand\bidu{^{\vee\vee}}
\newcommand{\leer}{\operatorname{--}}

\newcommand\catr{\underline{\operatorname{tr}}}

\newcommand\ap{^{?}}
\makeatletter
\def\namelabel#1#2{\@bsphack
  \protected@write\@auxout{}%
         {\string\newlabel{#1.nme}{{#2}{#2}}}%
  \@esphack}
\def\nmlabel#1#2{\label{#2}\namelabel{#2}{#1}}
\newcommand\nmref[1]{\ref{#1.nme}\ \ref{#1}}
\makeatother

\title{Higher Frobenius-Schur Indicators for Pivotal Categories}
\author{Siu-Hung Ng}
\address{Department of Mathematics, Iowa State University, Ames, IA 50011, USA}

\email{rng@iastate.edu}
\author{Peter Schauenburg}
\address{Mathematisches Institut der Universit\"at M\"unchen,
Theresienstr.\ 39, 80333 M\"unchen, Germany}
\email{schauenburg@math.lmu.de}

\begin{document}
\begin{abstract}
  We define higher Frobenius-Schur indicators for objects in
  linear pivotal monoidal categories. We prove that they are
  category invariants, and take values in the cyclotomic integers.
  We also define a family of natural endomorphisms of the identity
  endofunctor on a $k$-linear semisimple rigid monoidal category,
  which we call the Frobenius-Schur endomorphisms. For a
  $k$-linear semisimple pivotal monoidal category --- where both notions
  are defined ---, the
  Frobenius-Schur indicators can be computed as traces of
  the Frobenius-Schur endomorphisms.
\end{abstract}

\maketitle

\section*{Introduction}

The classical (degree two) Frobenius-Schur indicator $\nu_2(V)$ of
an irreducible representation $V$ of a finite group $G$ has been
generalized to Frobenius-Schur indicators of simple modules of
semisimple Hopf algebras by Linchenko and Montgomery
\cite{LinMon:FSTHA}, to certain $C^*$-fusion categories by Fuchs,
Ganchev, Szlach\'anyi, and Vescerny\'es \cite{FucGanSzlVes:FSIRMCC},
and further to simple objects in pivotal (or sovereign) categories
by Fuchs and Schweigert \cite{FucSch:CTCBC}; Mason and Ng
\cite{MasNg:CIFSISQHA} treated the case of simple modules over
semisimple quasi-Hopf algebras.
%and to simple modules of
%semisimple quasi-Hopf algebras by Mason and Ng
%\cite{MasNg:CIFSISQHA}.
As in the classical case, the indicator always takes one of the
values $0,\pm 1$, and is related to the question if and how the
representation in consideration is self-dual. In the last case,
proving that $\nu_2(V)\in\{0,\pm 1\}$ uses the result of Etingof,
Nikshych, and Ostrik \cite{EtiNikOst:FC} that the module category of
a semisimple complex quasi-Hopf algebra is a pivotal monoidal
category. A different proof based on this pivotal structure and the
description of indicators in \cite{FucGanSzlVes:FSIRMCC} was given
in \cite{Sch:FSIQHA}.

The higher indicators $\nu_n(V)$ of an irreducible group
representation $V$, which have less obvious meaning for the
structure of $V$, were generalized to simple modules of a semisimple
Hopf algebra by Kashina, Sommerh\"auser, and Zhu
\cite{KasSomZhu:HFSI}.

In the present paper we define and study higher Frobenius-Schur
indicators $\nu_n(V)$ for an object $V$ of a $k$-linear pivotal
monoidal category $\CC$. We do this with a view towards the
categories of modules over semisimple complex quasi-Hopf algebras,
though we will only give the (quite involved) explicit formulas
and examples for that case in another paper \cite{NS052}.

Our definition of $\nu_n(V)$ is as the trace of an endomorphism
$E_V^{(n)}$ of the vector space of morphisms $\CC(I,V^{\o n})$,
where $I$ is the unit object. The endomorphism arises from a special
case of a map $\CC(I,V\o W)\to\CC(I,W\o V)$ defined for any two
objects in terms of duality and the pivotal structure. In the case
where $V$ is an $H$-module for a semisimple Hopf algebra $H$, we can
identify $\CC(I,V^{\o n})$ with the invariant subspace $\left(V^{\o
n}\right)^H$, and $E_V^{(n)}$ is given by a cyclic permutation; the
description of $\nu_n(V)$ as a trace in this case is contained
in \cite{KasSomZhu:HFSI}. If $n=2$ and $V$ is simple, then
$\CC(I,V\o V)$ is one-dimensional or vanishes. Thus $E_V^{(2)}$ is
(at most) a scalar, which coincides with its trace; that computing
the trace in a different way leads to the indicator formulas from
\cite{MasNg:CIFSISQHA} was shown in \cite{Sch:FSIQHA}. An
endomorphism of $\CC(V\du,V)$ conjugate to $E_V^{(2)}$ is also used
to describe the degree two indicator in \cite{FucGanSzlVes:FSIRMCC}
(and to define $E_V^{(2)}$ in \cite{Sch:FSIQHA}). The maps
$E_V^{(n)}$ have been studied in connection with $3$-manifold
invariants by Gelfand and Kazhdan \cite{GelKaz:ITDM}.

We prove that the higher indicators are invariant under
equivalences of pivotal monoidal categories, and that equivalences
of pseudo-unital fusion categories (which are pivotal categories
by \cite{EtiNikOst:FC}) are automatically pivotal equivalences.
While these invariance properties are to be expected from the
categorical nature of the definitions, the proofs are not obvious.
In particular, it would be tautological to assume to be dealing
with strict categories (on grounds that every monoidal category is
\emph{equivalent} to a strict one) to prove invariance of certain
properties under \emph{equivalence}. Barrett and Westbury
\cite{BarWes:SC} have proved (although starting from a different
set of axioms) that every pivotal monoidal category is equivalent
to a strict one. By the invariance results we have proved, it is
then sufficient to assume such a simplified structure of $\CC$
when proving general properties of $E_V^{(n)}$ and $\nu_n(V)$. We
prove, in section 5, that the order of $E_V^{(n)}$ divides $n$;
this is stated in \cite{GelKaz:ITDM} and proved for $n=2$ (and
$\CC$ strict). In particular, the possible values of the higher
Frobenius-Schur indicators are cyclotomic integers; this is
well-known for the group case, where the higher indicators are in
fact always integers, and proved for the Hopf algebra case in
\cite{KasSomZhu:HFSI}. Finally, we show that the Frobenius-Schur
indicator of an object $V$ in a pivotal fusion category can be
computed as the trace (in a suitable sense) of a natural
endomorphism of $V$, which we call the Frobenius-Schur
endomorphism.

The organization of the paper is as follows: we cover in Section 1
some basic definitions, notations, conventions and preliminary
results of pivotal monoidal categories for the remaining discussion.
In Section 2, we give a proof that every pivotal monoidal category
is equivalent, as pivotal monoidal categories, to a {\em strict}
one. We then define a sequence of scalars $\nu_n(V)$, called the
{\em higher Frobenius-Schur indicators}, for each object $V$ in a
$k$-linear pivotal monoidal category in Section 3. We prove in
Section 4 that these higher Frobenius-Schur indicators are invariant
under $k$-linear pivotal monoidal equivalences. This result allows
one to study these indicators by considering {\em only} the strict
$k$-linear pivotal monoidal categories and we prove in Section 5
that all the higher Frobenius-Schur indicators for a $k$-linear
pivotal monoidal category are cyclotomic integers provided the
characteristic of the algebraically closed field $k$ is zero. In
Section 6, we consider the higher Frobenius-Schur indicators for a
{\em pseudo-unitary fusion} category over $\BC$. In this case, we
show that these indicators are \leaveout{indeed }invariant under
$k$-linearly monoidal equivalence. Finally, in Section 7, we show
that the $n$th Frobenius-Schur indicator of an object $V$ in a
semisimple $k$-linear pivotal monoidal category is\leaveout{ indeed
}the pivotal trace of a natural endomorphism $\FS_V^{(n)}$, called
the {\em Frobenius-Schur endomorphism}. In another paper
\cite{NS052} we will study the pivotal fusion categories of modules
over a semisimple complex quasi-Hopf algebra. In this case the
Frobenius-Schur endomorphisms correspond to central gauge invariants
in the quasi-Hopf algebras. The indicators are obtained by applying
the representation's character, which corresponds precisely to the
definition of higher indicators in \cite{KasSomZhu:HFSI} for the
case of ordinary Hopf algebras.

\subsection*{Acknowledgements}
The first author acknowledges support from the NSA grant number
H98230-05-1-0020. The second author thanks the \emph{Deutsche
Forschungsgemeinschaft} for support by a Heisenberg fellowship, the
Institute of Mathematics of Tsukuba University for its hospitality,
and Akira Masuoka for being the perfect host.

\section{Preliminaries}

We will first fix some conventions: In a monoidal category $\CC$,
the associativity  isomorphism is $\Phi\colon (U\o V)\o W\to U\o(V\o
W)$. We will make the assumption that the unit object in a monoidal
category is always strict, $V\o I=V=I\o V$. As pointed out in
\cite{Sch:TMCISO}, this assumption can always be made true after
replacing the tensor product in $\CC$ with an isomorphic one. For
the purpose of this paper however (where frequently the question of
invariance of certain constructions under tensor equivalences is
key) this may not be a rigorous justification; we simply make the
assumption that $I$ is a strict unit for simplicity, and hold that
the general case can be treated by an insignificant but annoying
expansion of all proofs. By the well-known coherence theorem for
monoidal categories, if $X,Y\in\CC$ are formed by tensoring the same
sequence of objects $V_1,\dots,V_n\in\CC$, only with different
placement of parentheses, then there is a unique morphism
$\Phi^{?}\colon X\to Y$ composed formally from instances of $\Phi$
and $\Phi\inv$. (By a ``formal'' composition of instances of $\Phi$
and $\Phi\inv$ we mean one that could be written down in a suitably
defined free category, excluding compositions that only become
possible because ``formally different'' objects happen to be
identical in the concrete category at hand; such accidental
composites may of course fail to agree.) As a simple example

\[\Phi\ap=\Phi\inv(T\o\Phi)=\Phi(\Phi\inv\o W)\Phi\inv\colon T\o((U\o V)\o W)\to (T\o U)\o (V\o W).\]

A monoidal functor $(\FF,\xi)\colon\CC\to\DD$ will have the
structure isomorphism $\xi\colon\FF(V)\o\FF(W)\to\FF(V\o W)$. We
will assume that the structure isomorphism for the unit objects is
the identity $\FF(I)=I$. Given a general monoidal functor
$(\FF,\xi,\xi_0)$ with $\xi_0\colon I\to\FF(I)$ not the identity,
this can always be achieved by replacing $\FF$ with an isomorphic
functor $\FF'$ whose object map is given by $\FF'(X)=\FF(X)$ if
$X\neq I$, and $\FF'(I)=I$. Let $X$ be obtained from tensoring a
sequence $V_1,\dots,V_n$ of objects of $\CC$ with some choice of
parentheses, and let $X'$ be obtained from $V_i'=\FF(V_i)$ in the
same way. We will use the following special case of coherence of
monoidal functors: There is a unique formal composition
$\xi\ap\colon X'\to \FF(X)$ of instances of $\xi$, and if $Y,Y'$
are obtained from the same sequence of objects with a different
placement of parentheses, then
\[\xymatrix{X'\ar[d]^{\Phi\ap}\ar[r]^-{\xi\ap}&\FF(X)\ar[d]^{\FF(\Phi\ap)}\\Y'\ar[r]^-{\xi\ap}&\FF(Y)}\]
commutes. The coherence result of Epstein \cite{Eps:FBTC} treats
the case of symmetric monoidal functors between symmetric
 monoidal categories; the diagram
above is contained in the appropriate analog for the non-symmetric
 case, obtained essentially by
just leaving out the symmetry. This is certainly folklore; a proof
can be extracted from \cite{Eps:FBTC}.

A (left) dual object of $V\in\CC$ is a triple $(V\du,\ev,\db)$
with an object $V\du\in\CC$ and morphisms $\ev\colon V\du\o V\to
I$ and $\db\colon I\to V\o V\du$ such that the compositions
\begin{gather*}
V\xrightarrow{\db\o V}(V\o V\du)\o V\xrightarrow{\Phi}V\o(V\du\o
V)\xrightarrow{V\o\ev}V,
\\
V\du\xrightarrow{V\du\o\db}V\du\o(V\o
V\du)\xrightarrow{\Phi\inv}(V\du\o V)\o V\du\xrightarrow{\ev\o
V\du}V\du
\end{gather*}
are identities. We say that $\CC$ is (left) rigid if every object
has a dual. We use the notation $\redu V$ for the symmetric notion
of a right dual of an object $V\in\CC$; this is the same as a left
dual in the category $\CC^\sym$ in which the tensor product is
defined in the reverse order.

For any object $V$ of a monoidal category having a dual object
$\dual V$, we obtain an adjunction $$A_0\colon \CC(U,V\o
W)\cong\CC(\dual V\o U,W)$$ by
\begin{gather*}
  A_0(f)=\left(\dual V\o U\xrightarrow{\dual V\o f}\dual V\o (V\o
  W)\xrightarrow{\Phi\inv}(\dual V\o V)\o W\xrightarrow{\ev\o
  W}W\right),\\
  A_0\inv(g)=\left(U\xrightarrow{\db\o U}(V\o\dual V)\o
  U\xrightarrow\Phi V\o(\dual V\o U)\xrightarrow{V\o g}V\o
  W\right).
\end{gather*}
In the case $U=I$ we abbreviate
\begin{equation}\label{ADef}
A_{V,W}=A_0 \colon \CC(I, V\o W) \lra \CC(\dual{V}, W).
\end{equation}
We will simply write $A$ for $A_{V, W}$ when the
context is clear.

 By the uniqueness of adjoints, dual objects in a monoidal
category are unique: If $(V\du,\ev,\db)$ and $(V^\circ,\ev',\db')$
are dual objects of $V$, then an isomorphism $v\colon V^\circ\to
V\du$ is uniquely determined by either of the conditions $\ev(V\o
v)=\ev'$ or $(V\o v)\db'=\db$.  An instance of this is the
compatibility of duals and tensor products: Given duals of $V$ and
$W$ in $\CC$, we get a dual $(W\du\o V\du,\tilde\ev,\tilde\db)$ of
$V\o W$, where
\begin{gather*}
\tilde\ev=\left((W\du\o V\du)\o (V\o
W)\xrightarrow{\Phi\ap}W\du\o((V\du\o V)\o
W)\xrightarrow{W\du\o\ev\o W}W\du\o W\xrightarrow{\ev}I\right),\\
\tilde\db=\left(I\xrightarrow{\db}V\o V\du\xrightarrow{V\o\db\o
V\du}V\o((W\o W\du)\o V\du)\xrightarrow{\Phi\ap}(V\o W)\o(W\du\o
V\du)\right).
\end{gather*}
If $(V\o W)\du$ is another dual of $V\o W$, this yields a unique
isomorphism $\zeta\colon W\du\o V\du\to (V\o W)\du$ with
$\tilde\ev=\ev_{V\o W}(\zeta\o V\o W)$ and $\db_{V\o W}=(V\o
W\o\zeta)\tilde\db$.

If $\CC$ is rigid, and we choose a dual object for each object in
$\CC$, then $(\leer)\du$ is a contravariant functor in a natural
way; we will always choose $I\du=I$, with both morphisms $\ev$ and
$\db$ identities. In this case the morphisms $\zeta$ are the
components of a monoidal functor structure on $(\leer)\du$.

Let $\mathcal C,\mathcal D$ be two rigid monoidal categories. Fix
a left dual $(X\du ,\ev,\db)$ for each object $X$ in $\mathcal C$
or $\mathcal D$, and make $(\leer)\du$ contravariant functors in
the natural way. Let $(\mathcal F,\xi)\colon \mathcal C\to\mathcal
D$ be a monoidal functor, with monoidal functor structure
$\xi\colon  \mathcal F(X)\o\mathcal F(Y)\to\mathcal F(X\o Y) $.
Then for each $X\in \CC$, we have a dual $\FF(X)^\circ=(\mathcal
F(X\du),\ev',\db')$ of $\mathcal F(X)$ in $\mathcal D$, with
\begin{gather*}\ev'=\left(\mathcal F(X\du)\o\mathcal F(X)\xrightarrow\xi\mathcal
F(X\du\o X)\xrightarrow{\mathcal
F(\ev)}I\right),\\
\db'=\left(I\xrightarrow{\FF(\db)}\FF(X\o
X\du)\xrightarrow{\xi\inv}\FF(X)\o\FF(X\du)\right).
\end{gather*}
We denote the canonical isomorphisms by
$\tilde\xi\colon\FF(X\du)\to\FF(X)\du$; they are the components of
a natural transformation which we call the \textbf{duality
transformation} of $(\FF,\xi)$.
\begin{lem}\nmlabel{Lemma}{dualtransmon}
  Let $(\FF,\xi)\colon\CC\to\DD$ be a monoidal functor between
  rigid monoidal categories. Then the duality transform
  $\tilde\xi$ is a monoidal natural transformation.
\end{lem}
\begin{proof}
  We want to show that
  $$\xymatrix{
    \FF(Y\du)\o\FF(X\du)
      \ar[r]^\xi
      \ar[d]^{\tilde\xi\o\tilde\xi}
    &\FF(Y\du\o X\du)
      \ar[r]^{\FF(\zeta)}
    &\FF((X\o Y)\du)
      \ar[d]^{\tilde\xi}
    \\\FF(Y)\du\o\FF(X)\du
      \ar[r]^{\zeta}
    &(\FF(X)\o\FF(Y))\du
      \ar[r]^-{(\xi\inv)\du}
    &\FF(X\o Y)\du}$$
  commutes. To do this, we compare the two morphisms
  $\FF(Y\du)\o\FF(X\du)\to (\FF(X)\o\FF(Y))\du$ arising from the
  diagram; note that two arrows $f,g\colon A\to B\du$ are the same
  iff $\ev_B(f \o B)=\ev_B(g \o B)$.

  On one hand we have the commutative diagram
  {\tiny\[\xymatrix{
    (\FF(Y\du)\o\FF(X\du))\o(\FF(X)\o\FF(Y))
      \ar[r]^-{\xi\o\id}
      \ar[d]_{\xi\ap}
    &\FF(Y\du\o X\du)\o(\FF(X)\o\FF(Y))
      \ar[r]^-{\FF(\zeta)\o\id}
    &\FF((X\o Y)\du)\o(\FF(X)\o\FF(Y))
      \ar[d]^{\xi\du\tilde\xi\o\id}
      \ar[dl]^{\xi\ap}
    \\\FF((Y\du\o X\du)\o(X\o Y))
      \ar[r]^-{\FF(\zeta\o\id)}
      \ar[dr]^{\FF(\tilde\ev)}
    &\FF((X\o Y)\du\o(X\o Y))
      \ar[d]^{\FF(\ev)}
    &(\FF(X)\o\FF(Y))\du\o(\FF(X)\o\FF(Y))
      \ar[dl]^{\ev}
    \\&I}\]}
  which, besides definitions and naturality, contains the diagram
  {\tiny\[\xymatrix{
    \FF((X\o Y)\du)\o(\FF(X)\o\FF(Y))
      \ar[r]^-{\tilde\xi\o\id}
      \ar[d]_{\id\o\xi}
    &\FF(X\o Y)\du\o(\FF(X)\o\FF(Y))
      \ar[d]_{\id\o\xi}
      \ar[r]^-{\xi\du\o\id}
    &(\FF(X)\o\FF(Y))\du\o(\FF(X)\o\FF(Y))
      \ar[dd]^{\ev}
    \\\FF((X\o Y)\du)\o\FF(X\o Y)
      \ar[r]^-{\tilde\xi\o\id}
      \ar[d]_\xi
    &\FF(X\o Y)\du\o\FF(X\o Y)
      \ar[dr]^\ev
    \\\FF((X\o Y)\du\o(X\o Y))
      \ar[rr]^-{\FF(\ev)}
    &&I}\]}
  in its lower right part. On the other hand the diagram
  \[\xymatrix{
    (\FF(Y\du)\o\FF(X\du))\o(\FF(X)\o\FF(Y))
      \ar[rr]^{\tilde\xi\o\tilde\xi \o \id}
      \ar[dr]^{\tilde\ev'}
      \ar[d]^{\xi\ap}
    &&(\FF(Y)\du\o\FF(X)\du)\o(\FF(X)\o\FF(Y))
      \ar[dl]_{\tilde\ev}
      \ar[d]^{\zeta\o\id}
    \\\FF((Y\du\o X\du)\o(X\o Y))
      \ar[r]_-{\FF(\tilde\ev)}
    &I
    &(\FF(X)\o\FF(Y))\du\o(\FF(X)\o\FF(Y)),
      \ar[l]^-{\ev}}\]
  in which $\tilde\ev'$ is analogous to $\tilde\ev$ defined above,
  but with two instances of $\ev'$ instead of $\ev$,
  commutes: Its upper triangle uses the definition of $\tilde\xi$
  twice, and its lower left triangle is an application of the
  definition of $\ev'$ and coherence; the diagram
  {\tiny\[\xymatrix{
    (\FF(Y\du)\o\FF(X\du))\o(\FF(X)\o\FF(Y))
      \ar[rr]^-{\Phi\ap}
      \ar[d]_{\xi\ap}
    &&\FF(Y\du)\o((\FF(X\du)\o\FF(X))\o\FF(Y))
      \ar[dl]_{\FF(Y\du)\o\xi\o\FF(Y)}
      \ar[d]^{\FF(Y\du)\o\ev'\o\FF(Y)}
    \\\FF((Y\du\o X\du)\o(X\o Y))
      \ar[d]_{\FF(\Phi\ap)}
    &\FF(Y\du)\o(\FF(X\du\o X)\o\FF(Y))
      \ar[r]^-{\FF(Y\du)\o\FF(\ev)\o\FF(Y)}
      \ar[dl]_{\xi\ap}
    &\FF(Y\du)\o\FF(Y)
      \ar[d]^{\ev'}
      \ar[dl]_\xi
    \\\FF(Y\du\o((X\du\o X)\o Y))
      \ar[r]^-{\FF(Y\du\o\ev\o Y)}
    &\FF(Y\du\o Y)
      \ar[r]^{\FF(\ev)}
    &I
    }\]}
  provides more details. Thus, we are done.
\end{proof}

A pivotal monoidal category is a rigid monoidal category with a
pivotal structure, that is, an isomorphism $j\colon V\to V\bidu$ of
monoidal functors. It follows that $j_{V\du}=\left(j_V\right)\inv$
for all $V\in\CC$, see \cite[Appendix]{Sch:FSIQHA}.

If $\CC,\DD$ are two pivotal monoidal categories, and $(\mathcal
F,\xi)\colon\CC\to\DD$ is a monoidal functor, we shall say that
$\mathcal F$ \textbf{preserves the pivotal structure}, if the
diagrams
$$\xymatrix{\FF(X)\ar[r]^-{\FF(j)}\ar[d]_j&\FF(\bidual
X)\ar[d]^{\tilde\xi}\\\bidual{\FF(X)}\ar[r]_-{\dual{\tilde\xi}}&\FF(X\du)\du}$$
commute.

Finally let us fix a few conventions on $k$-linear categories over a
field $k$: First of all, we will assume that the morphism spaces of
a $k$-linear monoidal category are finite-dimensional. A simple
object in a $k$-linear category $\CC$ is an object $V$ with
$\End_\CC(V)=k$ (note that sometimes such objects are called
absolutely simple in the literature). A semisimple $k$-linear
category is an abelian $k$-linear category in which every object is
a direct sum of simple objects. A semisimple $k$-linear monoidal
category is a monoidal category, and a semisimple $k$-linear
category, such that the tensor product is bilinear, and we will make
the general assumption that the unit object is simple.

It was noted in \cite{Ost:MCWHAMI}, see also
\cite[2.1]{EtiNikOst:FC}, that a semisimple $k$-linear left rigid
monoidal category is also right rigid. In fact if $\CC$ is
semisimple, then $\CC(I,V\o V\du)\neq 0$ implies $0\neq \CC(V\o
V\du,I)\cong\CC(V,V\bidu)$; thus any simple $V$ is isomorphic to
its bidual $V\bidu$, and $V\du$ is a right dual of $V$.

A fusion category \cite{BakKir:LTCMF,EtiNikOst:FC} is a semisimple
$\mathbb C$-linear rigid monoidal category with only finitely many
isomorphism classes of simple objects.

\section{Strictifying pivotal categories}

The main result in this section says that every pivotal category
can be assumed to be strict, which is supposed to mean that not
only its monoidal structure is strict, but also the isomorphisms
governing the compatibility between tensor product and duality as
well as the pivotal structure itself are identities. This is
proved by Barrett and Westbury in \cite{BarWes:SC}, who start from
a different set of axioms. The proof we give seems to be somewhat
shorter.

A pivotal strict monoidal category is a strict monoidal category
with a pivotal structure, that is, a rigid strict monoidal category
with a chosen fixed duality functor and a monoidal natural
isomorphism $j\colon X\to X\bidu$. By Mac Lane's coherence theorem,
in the form found, say, in Kassel's book \cite{Kas:QG}, every
pivotal monoidal category is equivalent, as a pivotal monoidal
category, to a pivotal strict monoidal category (one simply takes an
equivalent strict monoidal category and transports the pivotal
structure along the equivalence).

\begin{defn}
A \textbf{strict pivotal monoidal category} is a pivotal strict
monoidal category in which both the monoidal functor structure
$\zeta$ of $(\leer)\du$ (hence also that of $(\leer)\bidu$) and
the pivotal structure $j$ are identities.
\end{defn}

\begin{thm}\nmlabel{Theorem}{strict}
Every pivotal monoidal category is equivalent, as a pivotal
monoidal category, to a strict pivotal monoidal category.
\end{thm}
\begin{proof}
  Without loss of generality let $\CC$ be a pivotal strict monoidal
  category. We construct a strict pivotal monoidal category
  $\hat\CC$ and a pivotal strict monoidal equivalence $\FF\colon\hat\CC\to\CC$
  as follows:

  Objects of $\hat\CC$ are pairs $(\mathbf X,\mathbf\epsilon)$ in
  which $r\in{\mathbb N}_0$, $\mathbf X=(X_1,\dots,X_r)\in\CC^r$, and
  $\mathbf\epsilon=(\epsilon_1,\dots,\epsilon_r)\in\left({\mathbb Z}/(2)\right)^r$.
  The object map of $\FF$ is defined by $\FF(\mathbf X,\mathbf
  \epsilon)=X_1^{\vee\epsilon_1}\o\dots\o X_r^{\vee\epsilon_r}$,
  where we put $X^{\vee 0}=X$ and $X^{\vee 1}=X\du$ (and let the empty tensor product
  be the neutral object). We define
  \[\hat\CC((\mathbf X,\mathbf\epsilon),(\mathbf Y,\mathbf\delta))
  =\CC(\FF(\mathbf X,\mathbf\epsilon),\FF(\mathbf Y,\mathbf\delta))\]
  and let $\FF$ act on morphisms as the identity. Clearly $\FF$ is
  an equivalence, with a possible quasiinverse mapping $X\in\CC$ to
  $(X,1)\in\hat\CC$.

  We define the tensor product on $\hat\CC$ as componentwise
  concatenation on objects. With this choice,
  $\FF((\mathbf X,\mathbf\epsilon)\o(\mathbf Y,\mathbf\delta))
  =\FF(\mathbf X,\mathbf\epsilon)\o\FF(\mathbf Y,\mathbf\delta)$.
  In particular, we can define
  the tensor product of morphisms by
  taking the tensor product of morphisms in $\CC$; in
  this way $\hat\CC$ is a strict monoidal category, and $\FF$ is a
  strict monoidal functor.

  We choose a duality functor for $\hat\CC$ by putting
  \[(\mathbf X,\mathbf\epsilon)\du:=((X_r,\dots,X_1),(\epsilon_r+1,\dots,\epsilon_1+1))\]
  and defining evaluation inductively: for $r=0$, $\ev=\id_I$.
  For $\mathbb X\in\hat\CC$, $X\in\CC$ and $\epsilon\in{\mathbb Z/(2)}$ we define
  evaluation on $\mathbb X\o(X,\epsilon)$ by
  $$\ev=\left((X,\epsilon+1)\o\mathbb X\du\o\mathbb X\o(X,\epsilon)
  \xrightarrow{\id\o\ev\o\id}
  (X,\epsilon+1)\o(X,\epsilon)
  \xrightarrow{\ev}I\right),$$
  where $\ev\colon (X,\epsilon+1)\o(X,\epsilon)\to I$ is defined
  as
  \begin{align*}
  \left(X\du\o X\xrightarrow\ev I\right)&&\text{if }\epsilon=0,\\
  \left(X\o X\du\xrightarrow{j\o X\du}X\bidu\o
  X\du\xrightarrow{\ev}I\right)&&\text{if }\epsilon=1.
  \end{align*}
  It is clear that with a suitable choice of $\db$ this does define a dual object. Also, it
  follows from the definition that for any $\mathbb X,\mathbb
  Y\in\hat\CC$ we have
  \[\ev_{\mathbb X\o\mathbb Y}=\left(\mathbb Y\du\o\mathbb X\du\o \mathbb X\o\mathbb Y
    \xrightarrow{\id\o\ev\o\id}
    \mathbb Y\du\o\mathbb Y\xrightarrow{\ev}I\right),\]
  so that in $\hat\CC$ the duality functor is strict monoidal.
  Another direct consequence of the definition is that the
  component
  $$\tilde\xi_{(X,\epsilon)}\colon \FF(X,\epsilon+1)\to \FF(X,\epsilon)\du$$
  of the duality transformation $\tilde\xi$ associated to the
  trivial monoidal functor structure $\xi$ of the strict monoidal
  functor $\FF$ is given by
  \begin{align*}
    \id\colon X\du\to X\du&&\text{if }\epsilon=0,\\
    j\colon X\to X\bidu&&\text{if }\epsilon=1.
  \end{align*}

  We define a pivotal structure for $\hat\CC$ in the unique way
  that lets $\FF$ preserve the pivotal structure, that is, by
  requiring the diagrams
  $$\xymatrix{\FF(\mathbb X)\ar[r]^{\FF(j)}\ar[d]_j
    &\FF(\bidual{ \mathbb X})\ar[d]^{\tilde\xi}
    \\\bidual{\FF(\mathbb X)}\ar[r]_-{\dual{\tilde\xi}}
    &\FF(\mathbb X\du)\du}$$
  to commute for each $\mathbb X\in\hat\CC$. Observe that $\mathbb
  X\bidu=\mathbb X$ by definition. We wish to show that, moreover,
  the top arrow $\FF(j)$ is the identity for each $\mathbb X$.

  Since the pivotal structure is a monoidal transformation, it is enough
  to verify this for the
  case $\mathbb X=(X,\epsilon)$ with $X\in\CC$ and
  $\epsilon\in\mathbb Z/(2)$. In case $\epsilon=0$, the lower
  horizontal arrow is the identity, and the right vertical arrow
  is $j\colon X\to X\bidu$, so the top arrow is the identity.
  In case
  $\epsilon=1$, the right vertical arrow is the identity, and the
  lower horizontal arrow is $j\du$. But by \cite[Appendix]{Sch:FSIQHA} we have
  $j_X\du=j_{X\du}\inv$, so
  again the top arrow is the identity.
\end{proof}
\begin{remark}
  \nmref{dualtransmon} was used in the proof that the new pivotal structure
  is the identity, by reducing the question to the case of
  objects with only one component. Likely, the proof can be completed without using \nmref{dualtransmon}
  by an inductive argument.
\end{remark}

\section{Higher Frobenius-Schur indicators}
Throughout this section, $\CC$ is a pivotal monoidal category,
with pivotal structure $j$. We denote by $V^{\ro n}$ the $n$-fold
tensor power of an object $V\in\CC$ with rightmost parentheses;
thus $V^{\ro 0}=I$, and $V^{\ro(n+1)}=V\o V^{\ro n}$. There is a
unique isomorphism
$$\Phi^{(n)}\colon V^{\ro(n-1)}\o V\to V^{\ro n}$$ composed of
instances of $\Phi$; explicitly $\Phi^{(1)}$ is the identity, and
$$\Phi^{(n+1)}=\left(
  (V\o V^{\ro(n-1)})\o V
  \xrightarrow{\Phi}V\o(V^{\ro(n-1)}\o V)
  \xrightarrow{V\o\Phi^{(n)}}
  V^{\ro(n+1)}\right).$$
\begin{defn}
For $V,W\in\CC$, define $T_{VW}\colon\CC(\dual V,W)\to\CC(\dual
W,V)$ by
$$T_{VW}(f)=(\dual W\xrightarrow{\dual f}\bidual
V\xrightarrow{\piv_V\inv}V),$$ and put
\begin{gather*}E_{VW}=\left(\CC(I,V\o W)\xrightarrow{A}
  \CC(\dual V,W)\xrightarrow{T_{VW}}
  \CC(\dual W,V)\xrightarrow{A\inv}
  \CC(I,W\o V)\right),
  \\E^{(n)}_V=\left(\CC(I,V^{\ro n})\xrightarrow{E_{V,V^{\ro(n-1)}}}
                  \CC(I,V^{\ro{(n-1)}}\o
                  V)\xrightarrow{\CC(I,\Phi^{(n)})}
                  \CC(I,V^{\ro n})\right),
\end{gather*}
where $A$ denotes the adjunction \eqref{ADef} associated with
duality. Assume that $\CC$ is $k$-linear. Then for any positive
integers $r, n$, the $(n,r)$-th \textbf{Frobenius-Schur indicator}
of $V$ is the scalar
$$
\nu_{n,r}(V) = \Tr\left(\left(E_V^{(n)}\right)^r\right)\,.
$$
We will call $\nu_n(V):=\nu_{n,1}(V)$ the $n$-th
\textbf{Frobenius-Schur indicator} of $V$.
\end{defn}

We will give a reformulation of the definition of $E_{VW}$ which
will be useful later, in particular when we assume that the
categories in consideration are strict, and use graphical
notations. First, we define the following counterparts of the maps
$A_0$ and $A$:

\begin{defn}
We have an isomorphism $B_0\colon\CC(X\o V,W)\to\CC(X,W\o V\du)$
with
\begin{gather*}B_0(f)=\left(X\xrightarrow{X\o \db}X\o( V\o
V\du)\xrightarrow{\Phi\inv}(X\o V)\o V\du\xrightarrow{f\o V\du}W\o
V\du\right),\\
B_0\inv(g)=\left(X\o V\xrightarrow{g\o V}(W\o V\du)\o
V\xrightarrow{\Phi}W\o(V\du\o V)\xrightarrow{W\o\ev}W\right).
\end{gather*}
We note the special case $X=I$, with
$B=B_0\colon\CC(V,W)\to\CC(I,W\o V\du)$ given by
\begin{gather*}
B(f)=\left(I\xrightarrow\db V\o V\du\xrightarrow{f\o V\du}W\o
V\du\right),\\
B\inv(g)=\left(V\xrightarrow{g\o V}(W\o V\du)\o
V\xrightarrow{\Phi}W\o(V\du\o V)\xrightarrow{W\o\ev}W\right).
\end{gather*}
\end{defn}

It is easy to see (and well-known, at least implicit for example
in \cite{JoyStr:BTC}) that the composition
\[\CC(V,W)\xrightarrow{B}\CC(I,W\o V\du)\xrightarrow{A}\CC(W\du,V\du)\]
is given by $AB(f)=f\du$.

\begin{defn}\nmlabel{Definition}{DDef}
Denote by $D=D_{V,W}\colon \CC(I,V\o W)\to\CC(I,W\o V\bidu)$ the
composition
\[D=\left(\CC(I,V\o W)\xrightarrow{A}\CC(V\du,W)\xrightarrow{B}\CC(I,W\o V\bidu)\right).\]
\end{defn}
\begin{lem}\nmlabel{Lemma}{DLemma}We have
  \[E_{VW}=\left(\CC(I,V\o W)\xrightarrow{D_{VW}}\CC(I,W\o V\bidu)\xrightarrow{\CC(I,W\o j\inv)}\CC(I,W\o V)\right).\]
\end{lem}
\begin{proof}
  Comparing the definitions of $D$ and $E$, we have to show that
  the outer pentagon of the diagram
  \[\xymatrix{\CC(V\du,W)\ar[r]^-{B}\ar[d]_{(\leer)\du}
    &\CC(I,W\o V\bidu)\ar[dd]^{\CC(I,W\o j\inv)}\ar[dl]^{A}\\
    \CC(W\du,V\bidu)\ar[d]_{\CC(W\du,j\inv)}
    \\
    \CC(W\du,V)\ar[r]^{A\inv}&\CC(I,W\o V)}\]
  commutes. But we know that the triangle commutes, and the
  quadrangle is naturality of $A$ applied to $j\inv$.
\end{proof}

\section{Invariance I}
The definitions of the maps $E_{VW}$ and $E^{(n)}_V$ are given in
terms of monoidal structure, duality, and pivotal structure. Thus
if $\CC,\DD$ are two pivotal monoidal categories, and
$(\FF,\xi)\colon\CC\to\DD$ is monoidal (hence preserves duality),
and preserves the pivotal structure, then $\FF$ preserves $E_{VW}$
and $E^{(n)}_V$, in a sense we will make precise shortly. While
the statement is conceptually evident, the proof will involve some
technical machinery. The technical difficulties are due mostly to
the fact that $E_{VW}$ is a map of sets, and not a morphism
defined in the category.

\begin{lem}\nmlabel{Lemma}{FFA}
  Let $(\FF,\xi)\colon\CC\to\DD$ be a monoidal functor between
  rigid monoidal categories $\CC,\DD$. Then for all $V,W\in\CC$
  we have
  \[\xymatrix{
    \CC(I,V\o W)
      \ar[d]^-A
      \ar[r]^-\FF
    &\DD(I,\FF(V\o W))
      \ar[rr]^-{\DD(I,\xi\inv)}
    &&\DD(I,\FF(V)\o\FF(W))
      \ar[d]_A
    \\\CC(V\du,W)
      \ar[r]^-\FF
    &\DD(\FF(V\du),\FF(W))
      \ar[rr]^-{\DD(\tilde\xi\inv,\FF(W))}
    &&\DD(\FF(V)\du,\FF(W))}\]
  or $\FF(A(f))\tilde\xi\inv=A(\xi\inv\FF(f))$ for $f\colon I\to
  V\o W$.
\end{lem}
\begin{proof}
It is convenient to prove
$A\inv(\FF(g)\tilde\xi\inv)=\xi\inv\FF(A\inv(g))$ for $g\colon
V\du\to W$ instead, by a look at the diagram
\[\xymatrix{
  I\ar[d]_{\db}\ar[rr]^-{\FF(\db)}
  &&\FF(V\o V\du)\ar[rr]^-{\FF(V\o g)}\ar[d]_{\xi\inv}
  &&\FF(V\o W)\ar[d]_{\xi\inv}
  \\\FF(V)\o\FF(V)\du\ar[rr]^-{\FF(V)\o\tilde\xi\inv}
  &&\FF(V)\o\FF(V\du)\ar[rr]^-{\FF(V)\o\FF(g)}
  &&\FF(V)\o\FF(W)}\]
  whose top line is $\FF(A\inv(g))$.
\end{proof}

Similarly:
\begin{lem}\nmlabel{Lemma}{FFB}
  Let $(\FF,\xi)\colon\CC\to\DD$ be a monoidal functor between
  rigid monoidal categories $\CC,\DD$. Then for all $V,W\in\CC$
  we have
  \[\xymatrix{
    \CC(V,W)
      \ar[rr]^-\FF
      \ar[d]_B
    &&\DD(\FF(V),\FF(W))
      \ar[rr]^-B
    &&\DD(I,\FF(W)\o\FF(V)\du)
      \ar[d]_{\DD(I,\FF(W)\o\tilde\xi\inv)}
    \\\CC(I,W\o V\du)
      \ar[rr]^{\FF}
    &&\DD(I,\FF(W\o V\du))
      \ar[rr]^-{\DD(I,\xi\inv)}
    &&\DD(I,\FF(W)\o\FF(V\du))
      }\]
\end{lem}

Recall that there is a unique isomorphism
$\xi^{(n)}=\xi\ap\colon\FF(V^{\ro n})\to\FF(V)^{\ro n}$ composed
formally of instances of the monoidal functor structure $\xi$.
Explicitly $\xi^{0}=\id_I$, $\xi^{(1)}=\id_{\FF(V)}$, and
$$\xi^{(n+1)}=\left(
\FF(V^{\otimes (n+1)})\xrightarrow{\xi\inv_{V, V^{\otimes n}}}
\FF(V) \otimes \FF(V^{\otimes n}) \xrightarrow{V \otimes
\xi^{(n)}}\FF(V)^{\otimes (n+1)}\right)\,.
$$

\begin{prop}\nmlabel{Proposition}{p:E_invariance}
Let $(\FF,\xi)\colon\CC\to\DD$ be a  monoidal functor that
preserves the pivotal structure. Then for any $V, W \in \mathcal
C$, we have the following commutative diagrams:
\begin{equation}\label{EVW-inv}\begin{array}c
\xymatrix{
\CC(I, V \otimes W) \ar[rr]^-{E_{V, W}}\ar[d]^-{\FF} && \CC(I, W \otimes V)\ar[d]^-{\FF} \\
\DD(I, \FF(V \otimes W))\ar[d]_-{\DD(I, \xi\inv)}&& \DD(I, \FF(W
\otimes V))
\ar[d]^-{\DD(I, \xi\inv)}\\
\DD(I, \FF(V) \otimes \FF(W))\ar[rr]^-{E_{\FF(V), \FF(W)}} && \DD(I, \FF(W) \otimes \FF(V))\\
}\end{array}\end{equation}
\begin{equation}\label{EVn-inv}\begin{array}c
\xymatrix{\CC(I,V^{\ro n})\ar[rr]^-{E_V^{(n)}}\ar[d]_{\FF}
  &&\CC(I,V^{\ro n})\ar[d]^{\FF}\\
  \DD(I,\FF(V^{\ro n}))\ar[d]_-{\DD(I, \xi^{(n)})}
  &&\DD(I,\FF(V^{\ro n}))\ar[d]^-{\DD(I, \xi^{(n)})}\\
  \DD(I,\FF (V)^{\ro n})\ar[rr]_-{E_{\FF (V)}^{(n)}}
  &&\DD(I,\FF (V)^{\ro n})}
\end{array}\end{equation}
\end{prop}
\begin{proof}
  By \nmref{DLemma}, functoriality and naturality
  \[\xymatrix{
    &\CC(I,V\o W)
      \ar[dr]^E
      \ar[dl]_D
    \\
    \CC(I,W\o V\bidu)
      \ar[rr]^-{\CC(I,W\o j\inv)}
      \ar[d]_\FF
    &&\CC(I,W\o V)
      \ar[d]^\FF
    \\\DD(I,\FF(W\o V\bidu))
      \ar[rr]^-{\DD(I,\FF(W\o j\inv)}
      \ar[d]_{\DD(I,\xi\inv)}
    &&\DD(I,\FF(W\o V))
      \ar[d]^{\DD(I,\xi\inv)}
    \\\DD(I,\FF(W)\o\FF(V\bidu))
      \ar[rr]^-{\DD(I,\FF(W)\o\FF(j\inv))}
    &&\DD(I,\FF(W)\o\FF(V))
    }\]
  commutes. Thus, to prove \eqref{EVW-inv} we have to check commutativity of the right hand
  area of
  \[\xymatrix{
    \CC(V\du,W)
      \ar[r]^-{A\inv}
      \ar[d]_\FF
    &\CC(I,V\o W)
      \ar[r]^-D
      \ar[d]_\FF
    &\CC(I,W\o V\bidu)
      \ar[d]^\FF
    \\
    \DD(\FF(V\du),\FF(W))
      \ar[d]_{\DD(\tilde\xi\inv,\FF(W))}
    &\DD(I,\FF(V\o W))
      \ar[d]_{\DD(I,\xi\inv)}
    &\DD(I,\FF(W\o V\bidu))
      \ar[d]^{\DD(I,\xi\inv)}
    \\
    \DD(\FF(V)\du,\FF(W))
      \ar[r]^-{A\inv}
    &\DD(I,\FF(V)\o\FF(W))
      \ar[d]_D
    &\DD(I,\FF(W)\o\FF(V\bidu))
      \ar[d]^{\DD(I,\FF(W)\o\FF(j\inv))}
    \\
    &**[l]\DD(I,\FF(W)\o\FF(V)\bidu)
      \ar[r]^-{\DD(I,\FF(W)\o j\inv)}
    &\DD(I,\FF(W)\o\FF(V))}\]
  But the left hand area is \nmref{FFA}. The way around the
  outside of the diagram commutes by \nmref{FFB}, since $\FF$
  preserves the pivotal structure by assumption.

  \eqref{EVn-inv} follows from \eqref{EVW-inv} and the coherence of the monoidal functor
  $\FF$.
\end{proof}

The following corollary is an immediate consequence of Proposition
\ref{p:E_invariance}.
\begin{cor}\nmlabel{Corollary}{c:FS_invariance}
  Let $\CC$, $\DD$ be $k$-linear pivotal monoidal categories over a field $k$,
  and $\FF:\CC \lra \DD$ a $k$-linear monoidal equivalence that
  preserves the pivotal structure. Assume $\CC(I,V)$ and $\DD(I,W)$ are finite-dimensional
  for all objects $V$ in $\CC$ and $W$ in $\DD$. Then
  $$
  \nu_{n,r}(V)=\nu_{n,r}(\FF(V))\,.
  $$
  for any object $V$ in $\CC$ and positive integers $n, r$.
  \hfill\qed
\end{cor}

\section{The powers of $E$ and the values of $\nu_n$}

The main results of this section concern the powers of the map
$E_V^{(n)}$ --- whose $n$-th power turns out to be the identity
--- and the consequences for the possible values of its trace
$\nu_n$. In view of \nmref{p:E_invariance} and \nmref{strict}
we can assume that we are dealing with a strict pivotal monoidal
category to prove these results. Thus for the rest of this
section, except for the statements of the main results, we assume
that $\CC$ is a $k$-linear strict pivotal monoidal category. We
use graphical calculus in $\CC$ (see for example
\cite{JoyStr:GTCI} or \cite{Kas:QG}). In particular, we use
\[\gbeg33\got1{V\du}\gvac1\got1V\gnl\gwev3\gnl\gend,\qquad
\gbeg33\gnl\gwdb3\gnl\gob1 V\gvac1\gob1{V\du}\gend,\qquad\text{and
}\qquad\gbeg44\gnl\glmpb\gdnot f\gcmpb\gcmpb\grmpb\gnl
\gcl1\gcl1\gcl1\gcl1\gnl \gob1{V_1}\gob2{\cdots}\gob1{V_n}\gend\]
to depict the evaluation morphism $\ev\colon V\du\o V\to I$, the
dual basis morphism $\db\colon I\to V\o V\du$, and a morphism
$f\colon I\to V_1\o\dots\o V_n$, respectively.
 Thus, by \nmref{DDef} and \nmref{DLemma}
we have
\[
E_{VW}(f)=D_{VW}(f)=\gbeg45
        \gnl
        \gwdb4\gnl
        \gcl1\gwnot 2f\glmpb\grmpb\gcl2\gnl
        \gev\gcl1\gnl
        \gvac2\gob1W\gob1{V}
        \gend\]
for $f\colon I\to V\o W$. From this we deduce

\begin{equation}\label{Emonoidal}
  E_{V,W\o U}E_{U,V\o W}(f)
  =\gbeg65
    \gnl
    \gwdb6\gnl
    \gcl1\glmpb\gdnot{E_{U,V\o W}(f)}\gcmp\gcmpb\grmpb\gcl2\gnl
    \gev\gvac1\gcl1\gcl1\gnl
    \gvac3\gob1W\gob1U\gob1V\gend
  =\gbeg77
   \gnl
   \gwdb7\gnl
   \gcl3\gwdb5\gcl4\gnl
   \gvac1\gcl1\glmpb\gnot f\gcmpb\grmpb\gcl3\gnl
   \gvac1\gev\gcl1\gcl2\gcl2\gnl
   \gwev4\gnl
   \gvac4\gob1W\gob1U\gob1V\gend
  =E_{U\o V, W}(f)
\end{equation}
for each $f\colon I\to U\o V\o W$, which we use to prove
\begin{thm}
  Let $\CC$ be a pivotal monoidal category, $V\in\CC$ and $n\in\mathbb N$. Then
  $\left(E_V^{(n)}\right)^n=\id$.

  If $k$ is an algebraically closed field of
  characteristic zero and
  $\CC$ is a
  $k$-linear pivotal tensor category,
  then the  $(n, r)$-th Frobenius-Schur indicator
  $\nu_{n,r}(V)$
  is a cyclotomic integer in $\mathbb{Q}_n \subset k$ for any object $V \in \CC$
  and positive integers $n, r$, and we have
  $\nu_{n,n-r}(V)=\overline{\nu_{n,r}(V)}$.
\end{thm}
\begin{proof}
  It is sufficient to treat the case where $\CC$ is strict pivotal.
  We prove $\left(E_V^{(n)}\right)^k=E_{V^{\o k},V^{\o (n-k)}}$ for all $0\leq k\leq n$
  by induction:
  This is obvious for $k=0$ (and the definition of
  $E_V^{(n)}$ for $k=1$). Inductively,
  \begin{align*}
    \left(E_V^{(n)}\right)^{k+1}
      &=E_V^{(n)}\left(E_V^{(n)}\right)^k\\
      &=E_{V,V^{\o (n-1)}}E_{V^{\o k},V^{\o (n-k)}}\\
      &=E_{V,V^{\o (n-k-1)}\o V^{\o k}}E_{V^{\o k},V\o V^{\o
      (n-k-1)}}\\
      &=E_{V^{\o k}\o V,V^{\o (n-k-1)}}\\
      &=E_{V^{\o (k+1)},V^{\o (n-k-1)}}
  \end{align*}
  applying \eqref{Emonoidal} with $U=V^{\o k}$ and $W=V^{\o
  (n-k-1)}$.

  Now $\left(E_V^{(n)}\right)^n=E_{V^{\o n},I}$, so to show that
  $\left(E_V^{(n)}\right)^n=\id$, we are reduced to
  observing that both
  $A\colon\CC(I,V)\to\CC(V\du,I)$ and $B\colon\CC(V,I)\to\CC(I,V\du)$
  are special cases of the duality functor, and so
  $E_{V,I}(f)=f\bidu=f$ for all $f\colon I\to V$.

  If, in addition, $\CC$ is $k$-linear, then $E_V^{(n)}$ is a
  linear operator on the $k$-linear space $\CC(I, V^{\o n})$. Since
  $\left(E_V^{(n)}\right)^n=\id$, $\left(E_V^{(n)}\right)^r$ is diagonalizable for
  every
  integer $r$, and all its eigenvalues are $n$-th roots of unity. Thus,
  $$
  \nu_{n,r}(V)=\Tr\left(\left(E_V^{(n)}\right)^r\right)
  $$
  is a cyclotomic integer in $\mathbb{Q}_n$.

 Also, since $\left(E_V^{(n)}\right)^{n-r}$ is the inverse of
  $\left(E_V^{(n)}\right)^r$, the eigenvalues of these diagonalizable
  endomorphisms are inverse to each other. Since they are roots of
  unity, they are complex conjugate, and so are, consequently, the
  traces, which are the indicators $\nu_{n,n-r}(V)$, and
  $\nu_{n,r}(V)$, respectively.
\end{proof}
The definition of the map $E_V^{(n)}$ is not left-right symmetric;
we have chosen to ``bend the leftmost strand of $V^{\o n}$ over
the top''. Instead of repeating all the details of the definitions
in the opposite order, we will be very brief in giving the
relation between the two versions in the semisimple case:

  For $V\in\CC$ denote by $V^\sym$ the object $V$ viewed in the
  monoidal category $\CC^\sym$, and let $\nu_{n,k}({V^\sym})$ be
  the $(n,k)$-th Frobenius-Schur indicator of
  $V^\sym\in\CC^\sym$.

\begin{lem}\nmlabel{Lemma}{symmetric}
  Let $\CC$ be a semisimple $k$-linear pivotal monoidal category. Then
  $\nu_{n,k}({V^\sym})=\nu_{n,n-k}(V) $ for every $V\in\CC$ and
  integers $n,k$.
\end{lem}
\begin{proof}
  We can assume $1\leq k<n$.
  Since $\CC$ is semisimple, composition is a nondegenerate
  bilinear form $\CC(V^{\o n},I)\times\CC(I,V^{\o n})\to k$.
  Choose a basis $q_i\in\CC(I,V^{\o n})$ and let $p_i\colon V^{\o
  n}\to I$ be the elements of the dual basis. Thus the trace of an
  endomorphism $F\colon\CC(I,V^{\o n})\to\CC(I,V^{\o n})$ can be
  computed as $\Tr(F)=
  \sum_i p_iF(q_i)
  \in\CC(I,I)\cong k$.
  We use the graphical notations
  \[\gbeg64
    \got3{V^{\o(n-k)}}\got3{V^{\o k}}\gnl
    \gvac1\gcl1\gvac2\gcl1\gnl
    \gvac1\glmpt\gdnot {p}\gcmp\gcmp\grmpt\gnl
    \gend
    \qquad\text { and }\qquad
    \gbeg64
    \gnl
    \gvac1\glmpb\gdnot {q}\gcmp\gcmp\grmpb\gnl
    \gvac1\gcl1\gvac2\gcl1\gnl
    \gob3{V^{\o k}}\gob3{V^{\o(n-k)}}\gend\]
  for $q\in\CC(I,V^{\o n})$ and $p\in\CC(V^{\o n},I)$ to find
  \[\sum_i
  \gbeg46
  \gnl
  \gwdb4\gnl
  \gcl1\gdnot{q_i}\glmpb\grmpb\gcl2\gnl
  \gev\gcl1\gnl
  \gvac2\gdnot{p_i}\glmpt\grmpt\gnl
  \gnl\gend
  =
  \sum_i
  \gbeg68
  \gnl
  \gwdb6\gnl
  \gcl1\gwdb4\gcl4\gnl
  \gev\gdb\gcl1\gnl
  \gvac2\gcl2\gdnot{q_i\du}\glmpt\grmpt\gnl
  \gnl
  \gvac2\glmpt\gdnot{p_i}\gcmp\gcmp\grmpt\gnl
  \gend
  =
  \sum_i
  \gbeg68
  \gnl
  \gwdb6\gnl
  \gcl4\gwdb4\gcl1\gnl
  \gvac1\gcl1\gdb\gev\gnl
  \gvac1\gdnot{q_i\du}\glmpt\grmpt\gcl2\gnl
  \gnl
  \glmpt\gdnot{p_i}\gcmp\gcmp\grmpt\gnl
  \gend
  =
  \sum_i
  \gbeg46
  \gnl
  \gwdb4\gnl
  \gcl2\gdnot{q_i}\glmpb\grmpb\gcl1\gnl
  \gvac1\gcl1\gev\gnl
  \gdnot{p_i}\glmpt\grmpt\gnl
  \gend
\]
The left hand side computes the trace of
$\left(E_V^{(n)}\right)^k=E_{V^{\o k},V^{\o(n-k)}}$, while the
right hand side does the same for
$E_{(V^\sym)^{\o(n-k)},(V^\sym)^{\o
k}}=\left(E_{V^\sym}^{(n)}\right)^{n-k}$ in the category
$\CC^\sym$.
\end{proof}
\section{Invariance II}
Throughout this section we assume that $\CC$ is a fusion category.

The Grothendieck group $K(\CC)$ of $\CC$ is a ring with the
multiplication
$$
[U][V]=[U \o V]
$$
for any $U, V \in \CC$, where $[U]$ denotes the isomorphism class
of $U$. Note that $K(\CC)$ is a free abelian group with a basis
$\{[V_i]\}$ where $V_1, \dots, V_n$ is a complete set of
non-isomorphic simple objects of $\CC$. Thus, for any $V \in \CC$,
there exist unique non-negative integers $N_{i,1},\dots, N_{i, n}$
such that
$$
[V][V_i] = \sum_{j} N_{i,j}[V_j]\,.
$$
The Frobenius-Perron dimension of $V$, denoted by $\FPdim(V)$, is
defined to be the largest nonnegative real eigenvalue of the
non-negative integer matrix
$$
\rho_V = \left[N_{i,j}\right];
$$
it dominates the absolute values of all eigenvalues of $\rho_V$
(cf. \cite[Section 8]{EtiNikOst:FC}). The Frobenius-Perron
dimension $\FPdim(\CC)$ of $\CC$ is defined to be
$$
\sum_{i} \FPdim\left(V_i\right)^2\,.
$$
By \cite[Proposition 2.1]{EtiNikOst:FC}, for any simple object $V
\in \CC$, there exists an isomorphism $a: V \lra \bidual{V}$. One
can define the categorical trace $\catr(a) \in \mathbb C$ of $a$
by
$$
\catr(a)=(I \xrightarrow{\db} V \o \dual{V} \xrightarrow{a \o V}
\bidual{V}\o V\du \xrightarrow{\ev} I)\,.
$$
The normed square of $V$, denoted by $|V|^2$, is defined to be
$$
|V|^2 = \catr(a)\catr(\dual{(a\inv)})\,.
$$
Note that the definition of $|V|^2$ is independent of the choice
of $a$. The global dimension $\dim(\CC)$ of $\CC$ is defined to be
$$
\sum_i |V_i|^2\,.
$$
A fusion category $\CC$ over $\mathbb C$ is called pseudo-unitary
if $\dim(\CC)=\FPdim(\CC)$. It was shown in \cite[Proposition
8.23]{EtiNikOst:FC} that if $\CC$ is a pseudo-unitary fusion
category, there exists a unique pivotal structure $j:Id \lra
\bidual{(\leer)}$ such that
\begin{equation}\label{eq:canonical_pivotal_structure}
\catr(j_V)=\FPdim(V)
\end{equation}
for every simple object $V \in \CC$.

The following lemma says that monoidal functors preserve traces in
a suitable sense:
\begin{lem}
If $(\mathcal
  F,\xi)\colon\CC\to\DD$ is a monoidal functor between rigid $k$-linear monoidal categories,
  and $a\colon V\to
  V\bidu$ in $\CC$, then
  \begin{equation}
    \FF(\catr_\CC(a))=\catr_\DD\left(\FF(V)\xrightarrow
    {\FF(a)}\FF(V\bidu)\xrightarrow{\tilde\xi}\FF(V\du)\du\xrightarrow{(\tilde\xi\inv)\du}\FF(V)\bidu\right)
    .
  \end{equation}
\end{lem}
\begin{proof}
  We can write the categorical trace map as the composition
  \[\catr_\CC=\left(\CC(V,V\bidu)\xrightarrow{B_0\inv}\CC(V\o V\du,I)\xrightarrow{\CC(\db,I)}\CC(I,I)\right),\]
  so
  \[\catr_\CC(a)=\left(I\xrightarrow\db V\o V\du\xrightarrow{B_0\inv(a)}I\right).\]
  Writing $W=V\du$, the diagram
  \[\xymatrix{
    \FF(V)\o\FF(W)
      \ar[d]_{\xi}
      \ar[rr]^-{\FF(a)\o\FF(W)}
    &&\FF(W\du)\o\FF(W)
      \ar[rr]^-{\tilde\xi\o\FF(W)}
      \ar[d]^{\xi}
    &&\FF(W)\du\o\FF(W)
      \ar[d]_\ev
    \\\FF(V\o W)
      \ar[rr]^-{\FF(a\o W)}
    &&\FF(W\du\o W)
      \ar[rr]^-{\FF(\ev)}
    &&I
    }\]
  shows  $\FF(B_0\inv(a))\xi=B_0\inv(\tilde\xi\FF(a))$ for $a\colon
  V\to W\du$ in $\CC$, which is yet another variant of \nmref{FFA}.
  Further, $B_0\inv\colon\DD(X,Y\du)\to \DD(X\o
  Y,I)$ is natural, so for $X$ and $f\colon Y\to Z$ in $\DD$
  \[\begin{array}c
  \xymatrix{\DD(X,Z\du)\ar[r]^-{B_0\inv}\ar[d]_{\DD(X,f\du)}
    &\DD(X\o Z,I)\ar[d]^{\DD(X\o f,I)}
    \\\DD(X,Y\du)\ar[r]^-{B_0\inv}
    &\DD(X\o Y,I)}
  \end{array}
  \]
  commutes, that is to say, for each $b\colon X\to Z\du$, the diagram
  \[
  \begin{array}c
    \xymatrix{X\o Y\ar[r]^-{X\o f}\ar[dr]_{B_0\inv(f\du b)}&X\o Z\ar[d]^{B_0\inv(b)}\\&I}
  \end{array}
  \]
  commutes.
  As a consequence,
  \[\xymatrix{
    I
      \ar[rr]^-{\FF(\db)}
      \ar[dd]_\db
    &&\FF(V\o V\du)
      \ar[dd]^{\FF(B_0\inv(a))}
    \\&\FF(V)\o\FF(V\du)
      \ar[ur]^{\xi}
      \ar[dr]^{B_0\inv(\tilde\xi\FF(a))}
    \\\FF(V)\o\FF(V)\du
      \ar[ur]^{\FF(V)\o\tilde\xi\inv}
      \ar[rr]_-{B_0\inv((\tilde\xi\inv)\du\tilde\xi\FF(a))}
    &&I}\]
  commutes, and the outer rectangle is our claim.
\end{proof}
\begin{cor}
  Let $(\FF,\xi)\colon\CC\to\DD$ be a $\mathbb C$-linear monoidal equivalence of
  pseudo-unitary categories. Then $\FF$ preserves the canonical
  pivotal structures.
\end{cor}
\begin{proof}
  Being a $\mathbb C$-linear monoidal equivalence, $\FF$ preserves the
  Frobenius-Perron dimensions, so for any simple $V$ in $\CC$
  \[\catr_\DD((\tilde\xi\inv)\du\tilde\xi \FF(j_V))
  =\FF(\catr_\CC(j_V))
  =\catr_\CC(j_V)
  =\FPdim(V)
  =\FPdim(\FF(V))
  =\catr_\DD(j_{\FF(V)}).\]
  Since $\DD(W,W\bidu)$ is one-dimensional for any simple object $W$,
  the canonical pivotal structure on
  $\FF(V)$ is determined by its trace, and we are done.
\end{proof}

\begin{cor}
  The Frobenius-Schur indicators of the simple objects of a
  pseudo-unitary fusion category $\CC$ are invariants of $\CC$ as
  a monoidal category.
\end{cor}

\section{Frobenius-Schur endomorphisms}

Throughout this section $\CC$ is a semisimple $k$-linear monoidal
category.  Recall that in our conventions this includes the
condition that the unit object is simple, and that being simple
means that the endomorphism ring is isomorphic to the base field
$k$.

Let $S$ be a simple object in $\CC$. We say that an object $V$ is
$S$-isotypical if it is isomorphic to a sum of copies of $S$. An
$I$-isotypical object will be called a trivial object. The
subcategory of $\CC$ consisting of all trivial objects will be
denoted by $\Triv(\CC)$. An arbitrary object $V$ is isomorphic to
a direct sum of $S$-isotypical objects. More explicitly, there are
$S$-isotypical objects $V^{(S)}$ and morphisms $\iota^{(S)}\colon
V^{(S)}\to V$ and $\pi^{(S)}\colon V\to V^{(S)}$ for each simple
$S$ from a set of representatives of the isomorphism classes of
simples in $\CC$ such that
\begin{gather*}
\forall S\colon\qquad\pi^{(S)}\iota^{(S)}=\id_{V^{(S)}},\\
\sum_{S}\iota^{(S)}\pi^{(S)}=\id_V.
\end{gather*}
We will call the triple $(V^{(S)},\iota^{(S)},\pi^{(S)})$ an
$S$-isotypical component of $V$. Isotypical components are
uniquely determined up to isomorphism in an obvious sense. We will
only be needing $I$-isotypical components, which we will call
trivial components and denote by $(V^\triv,\iota,\pi)$.

\begin{lem}\nmlabel{Lemma}{trivialbraid}
  Let $\CC$ be a semisimple $k$-linear monoidal category.
  There is a unique
  isomorphism $\tau_{VT}\colon V\o T\to T\o V$ natural in
  $T\in\Triv(\CC)$ and $V\in\CC$ for which $\tau_{VI}$ is the
  identity (the canonical isomorphism).
  For all $V,W$ in $\CC$ the diagram
  \begin{equation}\label{braidingaxiom}
    \begin{array}c\xymatrix{(V\o W)\o T\ar[r]^-{\tau}\ar[d]_{\Phi}
    &T\o(V\o W)\ar[d]^{\Phi\inv}
    \\V\o(W\o T)\ar[d]_{V\o\tau}
    &(T\o V)\o W\ar[d]^{\tau\o W}
    \\V\o(T\o W)\ar[r]^-{\Phi\inv}
    &(V\o T)\o W}
    \end{array}\end{equation}
  commutes.
\end{lem}
\begin{proof}
  Fix $V$, and write $T$ as a direct sum of copies of $I$, with
  projections $\pi_i\colon T\to I$ and injections $\iota_i\colon
  I\to T$ for $i$ in some finite index set. If $\tau$ is to be a
  natural transformation, then the diagrams
  \[\xymatrix{V\o T\ar[r]^\tau\ar[d]_{V\o \pi_i}&T\o V\ar[d]^{\pi_i\o V}\\V\o I\ar@{=}[r]&I\o V}\]
  should commute, and there exists a unique isomorphism $\tau$ making
  these diagrams commute. Let $T'$ be another trivial object, and
  choose a direct sum decomposition with projections $\pi'_i\colon
  T'\to I$ and $\iota'_i\colon I\to T'$; construct $\tau_{VT'}$ analogous to $\tau_{TV}$.
  Let $f\colon T\to T'$ be
  a morphism, and let $\alpha_{ij}=\pi'_if\iota_j\colon I\to I$.
  The diagrams
  \[\xymatrix{V\o T\ar[r]^-{V\o f}
    &V\o T'\ar[r]^{\tau'}\ar[d]^{V\o \pi'_{i}}
    &T'\o V\ar[d]^{\pi'_i\o V}
    \\V\o I\ar[u]_{V\o \iota_j}\ar[r]_-{V\o \alpha_{ij}}
    &V\o I\ar@{=}[r]&I\o V}\]
  \[\xymatrix{V\o T\ar[r]^\tau
    &T\o V\ar[r]^{f\o V}
    &T'\o V\ar[d]^{\pi'_i\o V}
    \\V\o I\ar[u]_{V\o \iota_j}\ar@{=}[r]
    &I\o V\ar[u]_{\iota_j\o V}\ar[r]_{\alpha_{ij}\o V}
    &I\o V}\]
  and the fact that $\alpha_{ij}\o V=V\o\alpha_{ij}$ show that
  $(f\o V)\tau=\tau'(V\o f)$, so that $\tau$ is independent of the
  choice of the direct summand decompositions, and natural in $T$.
  Naturality in $V$ is obvious, and the braiding equality
  \eqref{braidingaxiom}
  follows from uniqueness.
\end{proof}
\begin{defn}
  Let $\CC$ be a semisimple $k$-linear rigid monoidal category, $V$ an object of $\CC$,
  and $n\in\mathbb N$.

  Choose a right dual $X$ of $V^{\o(n-1)}$, and a trivial
  component $((V^{\o n})^\triv,\iota,\pi)$.

  The $n$-th \textbf{Frobenius-Schur endomorphism} of $V$ is the
  composition
  \begin{multline*}\FS^{(n)}_V=\biggl(V\xrightarrow{\Unit}X\o(V^{\o(n-1)}\o V)
     \xrightarrow{X\o\Phi^{(n)}}X\o V^{\o n}
  \\     \xrightarrow{X\o\pi}X\o(V^{\o n})^\triv
     \xrightarrow{\tau}(V^{\o n})^\triv\o X
     \xrightarrow{\iota\o X}V^{\o n}\o X
     \xrightarrow{\Counit}V\biggr),
  \end{multline*}
  where
  \begin{gather*}
    \Unit=\left(V\xrightarrow{\db\o V}(X\o V^{\o(n-1)})\o V
      \xrightarrow{\Phi}X\o(V^{\o(n-1)}\o V)\right),\\
    \Counit=\left((V\o V^{\o(n-1)})\o
    X\xrightarrow{\Phi}V\o(V^{\o(n-1)}\o
    X)\xrightarrow{V\o\ev}V\right),
  \end{gather*}
  and $\tau$ is as in \nmref{trivialbraid}.
\end{defn}
\begin{lem}
  Let $V$ be an object of $\CC$.
  The $n$-th Frobenius-Schur endomorphism of $V$ is independent of
  the choice of a right dual of $V^{\o(n-1)}$, and of the choice
  of a trivial component functor. If $\FF\colon\CC\to\DD$ is an
  equivalence of semisimple $k$-linear monoidal categories, then
  $\FF(\FS^{(n)}_V)=\FS^{(n)}_{\FF(V)}$.
\end{lem}
\begin{proof}
By contrast to our earlier definitions of $E_{V}^{(n)}$ and the
indicators, the construction of $\FS^{(n)}$ takes place entirely
within the monoidal category $\CC$; the verifications are thus
routine, and we will omit the details.
\end{proof}
\begin{remark}
  Let $V$ be an object of $\CC$, and choose a
  trivial component $(V^\triv,\iota,\pi)$.
  Choose a direct sum decomposition of
  $V^\triv$ with projections $\pi_i\colon V^\triv\to S$ and
  injections $\iota_i\colon S\to V^\triv$ for $i$ from some finite
  index set.

  The composition $\CC(V,I)\times\CC(I,V)\to \CC(I,I)\cong k$ is a
  nondegenerate bilinear form, with $p_i=\pi_i\pi\colon V\to I$
  and $q_i=\iota\iota_i\colon I\to V$ as dual bases.
\end{remark}
The following definition follows Barrett and Westbury
\cite{BarWes:SC}:
\begin{defn}
Let $\CC$ be a pivotal monoidal category with pivotal structure
$j$. For an endomorphism $f\colon V\to V$ we define the left and
right \textbf{pivotal traces} to be
\begin{gather*}
\ptrr(f)=\catr(j_Vf)=\left(I\xrightarrow{\db}V\o
V\du\xrightarrow{f\o V\du}V\o V\du\xrightarrow{j_V\o V\du}V\bidu\o
V\du\xrightarrow{\ev}I\right)
\\
\ptrl(f)=\left(I\xrightarrow{\db}V\du\o V\bidu\xrightarrow{V\du\o
j_V\inv}V\du\o V\xrightarrow{V\du\o f}V\du\o
V\xrightarrow{\ev}I\right).
\end{gather*}
The category $\CC$ is called \textbf{spherical} if
$\ptrl(f)=\ptrr(f)=:\ptr(f)$ for all $f$.
\end{defn}
Note that $\ptrl(f)=\ptrr(f\du)$. It is well-known and easy to check
that the pivotal traces are cyclic, for $g\colon V\to W$ and
$f\colon W\to V$:
\begin{align*}
\ptrr(fg)&=\left(I\xrightarrow\db V\o V\du\xrightarrow{jfg\o
V\du}V\bidu\o V\du\xrightarrow\ev I\right)\\
&=\left(I\xrightarrow\db V\o V\du\xrightarrow{f\bidu jg\o
V\du}V\bidu\o V\du\xrightarrow\ev I\right)\\
&=\left(I\xrightarrow\db V\o V\du\xrightarrow{jg\o f\du}W\bidu\o
W\du\xrightarrow\ev I\right)\\
&=\left(I\xrightarrow\db W\o W\du\xrightarrow{jgf\o W\du}W\bidu\o
W\du\xrightarrow\ev I\right)\\
&=\ptrr(gf).
\end{align*}

In the case where $\CC$ is strict, the left and right pivotal
traces are given, in the graphical calculus, by closing the
morphism $f$ to a loop on the left, resp.\ on the right:
\[\ptrl(f)=\gbeg25\gnl
           \gdb\gnl
           \gcl1\gbmp f\gnl
           \gev\gnl
           \gend
 \qquad\text{ resp. }\qquad
 \ptrr(f)=\gbeg25\gnl
          \gdb\gnl
          \gbmp f\gcl1\gnl
          \gev\gnl
          \gend\]
\begin{thm}\nmlabel{Theorem}{traceformula}
  Let $\CC$ be a semisimple $k$-linear pivotal monoidal category. For an
  object $V$ in $\CC$ we have
  \[\nu_n(V)=\ptrl(\FS^{(n)}_V).\]
  In addition, if $\CC$ is spherical, then
 $\nu_n(V)=\nu_n(V\du)$.
\end{thm}
\begin{proof}
In view of the previous results, it is enough to prove the claim
under the assumption that $\CC$ is strict as a pivotal category.

Choose a trivial component $((V^{\o n})^\triv,\iota,\pi)$, and
choose a direct sum decomposition of $(V^{\o n})^\triv$ with
projections $\pi_i\colon (V^{\o n})^\triv\to I$ and injections
$\iota_i\colon I\to(V^{\o n})^\triv$. Put $p_i=\pi_i\pi$ and
$q_i=\iota\iota_i$.

Put $W=V^{\o(n-1)}$ and $T=(V^{\o n})^\triv$, and use the
graphical notations
\[\gbeg24
  \got1W\got1V\gnl
  \gcl1\gcl1\gnl
  \gdnot{p_i}\glmpt\grmpt\gnl
  \gend
  {\quad,}\qquad
  \gbeg24
  \gnl
  \gdnot{q_i}\glmpb\grmpb\gnl
  \gcl1\gcl1\gnl
  \gob1V\gob1W\gend
  {\quad,}\qquad
  \gbeg25
  \got1W\got1V\gnl
  \gcl1\gcl1\gnl
  \gdnot\pi\glmptb\grmpt\gnl
  \gcl1\gnl
  \gob1T\gend
  {\quad,}\qquad
  \gbeg25
  \gvac1\got1T\gnl
  \gvac1\gcl1\gnl
  \gdnot\iota\glmpb\grmptb\gnl
  \gcl1\gcl1\gnl
  \gob1V\gob1W\gend
  {\quad,}\qquad
  \tau_{WT}=
  \gbeg23
  \got1W\got1T\gnl
  \gcn1113\gcn111{-1}\gnl
  \gob1T\gob1W\gend\]
  to find
\begin{multline*}
\Tr(E_V^{(n)})
  =
  \sum_i
  \gbeg46
  \gnl
  \gwdb4\gnl
  \gcl1\gdnot{q_i}\glmpb\grmpb\gcl2\gnl
  \gev\gcl1\gnl
  \gvac2\gdnot{p_i}\glmpt\grmpt\gnl
  \gnl\gend
= \sum_i
  \gbeg58
  \gnl
  \gwdb5\gnl
  \gcl4\gvac1\gdb\gcl1\gnl
  \gvac2\gcl1\gdnot{p_i}\glmpt\grmpt\gnl
  \gvac2\gcn2113\gnl
  \gvac1\gdnot{q_i}\glmpb\grmpb\gcl1\gnl
  \gev\gev\gnl
  \gnl
  \gend
  =
  \gbeg58
  \gnl
  \gwdb5\gnl
  \gcl4\gvac1\gdb\gcl1\gnl
  \gvac2\gcl1\gdnot{\pi}\glmpt\grmpt\gnl
  \gvac2\gcn1113\gcn111{-1}\gnl
  \gvac1\gdnot{\iota}\glmpb\grmpb\gcl1\gnl
  \gev\gev\gnl
  \gnl
  \gend
  =
  \ptrl\left(
  \gbeg47
  \gvac3\got1V\gnl
  \gvac1\gdb\gcl1\gnl
  \gvac1\gcl1\gdnot{\pi}\glmptb\grmpt\gnl
  \gvac1\gcn1113\gcn111{-1}\gnl
  \gdnot\iota\glmpb\grmptb\gcl1\gnl
  \gcl1\gev\gnl
  \gob1V\gend\right){\quad.}
\end{multline*}
Note that
$$
(\FS_V^{(n)})\du =
   \gbeg79
  \got1{V\du}\gnl
  \gcl6\gvac3\gdb\gnl
  \gvac2\gdb\gcl1 \gcl6\gnl
  \gvac2\gcl1\gdnot{\pi}\glmpt\grmpt\gnl
  \gvac2\gcn1113\gcn111{-1}\gnl
  \gvac1\gdnot{\iota}\glmpb\grmpb\gcl1\gnl
   \gvac1\gcl1\gev\gnl
   \gev\gnl
   \gvac5\gob1{V\du} \gnl
  \gend
  =
  \gbeg57
  \gvac3\got1{V\du}\gnl
  \gvac1\gdb\gcl1\gnl
  \gvac1\gcl1\gdnot{\iota\du}\glmptb\grmpt\gnl
  \gvac1\gcn1113\gcn111{-1}\gnl
  \gdnot{\pi\du}\glmpb\grmptb\gcl1\gnl
  \gcl1\gev\gnl
  \gob1{V\du}\gend
  =\quad\FS_{V\du}^{(n)}\,.
$$
If $\CC$ is spherical, we have
$$
\nu_n(V)=\ptrl(\FS_V^{(n)})=\ptrl((\FS_V^{(n)})\du) =
\ptrl(\FS_{V\du}^{(n)})=\nu_n(V\du)\,. \qedhere
$$
\end{proof}
\begin{prop}
The Frobenius-Schur endomorphism is a natural endomorphism of the
identity functor on a semisimple $k$-linear rigid monoidal
category.
\end{prop}
\begin{proof}
It suffices to prove this in the case where $\CC$ is strict.
Consider $f\colon V\to W$ in $\CC$. Write
\[f_k=W^{\o(k-1)}\o f\o V^{\o(n-k)}\colon W^{\o(k-1)}\o V^{\o(n-k+1)}\to W^{\o k}\o V^{\o(n-k)},\]
and denote by $\iota_k,\pi_k$ the inclusion and projection of
$W^{\o k}\o V^{\o (n-k)}$ from and to its trivial component. In
particular $f_k\iota_{k-1}=\iota_kf_k^\triv$ and
$\pi_kf_k=f_k^\triv\pi_{k-1}$. Now for $1\leq k\leq n$
\[\gbeg 49
  \gvac3\got1V\gnl
  \gvac1\gdb\gcl1\gnl
  \gvac1\gcl1\gdnot{\pi_{k-1}}\glmptb\grmpt\gnl
  \gvac1\gsym\gnl
  \gdnot{\iota_{k-1}}\glmpb\grmptb\gcl3\gnl
  \gcl1\gcl1\gnl
  \gdnot{f_k}\glmptb\grmptb\gnl
  \gcl1\gev\gnl
  \gob1W\gend
  =
  \gbeg49
  \gvac3\got1V\gnl
  \gvac1\gdb\gcl1\gnl
  \gvac1\gcl1\gdnot{\pi_{k-1}}\glmptb\grmpt\gnl
  \gvac1\gsym\gnl
  \gdnot{f_k^\triv}\glmp\grmptb\gcl3\gnl
  \gvac1\gcl1\gnl
  \gdnot{\iota_k}\glmpb\grmptb\gnl
  \gcl1\gev\gnl
  \gob1W\gend
  =
  \gbeg49
  \gvac3\got1V\gnl
  \gvac1\gdb\gcl1\gnl
  \gvac1\gcl3\gdnot{\pi_{k-1}}\glmptb\grmpt\gnl
  \gvac2\gcl1\gnl
  \gvac2\gdnot{f_k^\triv}\glmptb\grmp\gnl
  \gvac1\gsym\gnl
  \gdnot{\iota_k}\glmpb\grmptb\gcl1\gnl
  \gcl1\gev\gnl
  \gob1W\gend
  =
  \gbeg49
  \gvac3\got1V\gnl
  \gvac1\gdb\gcl1\gnl
  \gvac1\gcl3\gdnot{f_k}\glmptb\grmptb\gnl
  \gvac2\gcl1\gcl1\gnl
  \gvac2\gdnot{\pi_k}\glmptb\grmpt\gnl
  \gvac1\gsym\gnl
  \gdnot{\iota_k}\glmpb\grmptb\gcl1\gnl
  \gcl1\gev\gnl
  \gob1W\gend
  \]
  and for $1\leq k<n$
  \[\gbeg49
  \gvac3\got1V\gnl
  \gvac1\gdb\gcl1\gnl
  \gvac1\gcl3\gdnot{f_k}\glmptb\grmptb\gnl
  \gvac2\gcl1\gcl1\gnl
  \gvac2\gdnot{\pi_k}\glmptb\grmpt\gnl
  \gvac1\gsym\gnl
  \gdnot{\iota_k}\glmpb\grmptb\gcl1\gnl
  \gcl1\gev\gnl
  \gob1W\gend
  =
  \gbeg49
  \gvac3\got1V\gnl
  \gvac1\gdb\gcl3\gnl
  \gvac1\gcl3\gbmp{f_k}\gnl
  \gvac2\gcl1\gnl
  \gvac2\gdnot{\pi_k}\glmptb\grmpt\gnl
  \gvac1\gsym\gnl
  \gdnot{\iota_k}\glmpb\grmptb\gcl1\gnl
  \gcl1\gev\gnl
  \gob1W\gend
  =
  \gbeg49
  \gvac3\got1V\gnl
  \gvac1\gdb\gcl3\gnl
  \gdnot{\redu{(f_k)}}\glmp\grmptb\gcl2\gnl
  \gvac1\gcl2\gnl
  \gvac2\gdnot{\pi_k}\glmptb\grmpt\gnl
  \gvac1\gsym\gnl
  \gdnot{\iota_k}\glmpb\grmptb\gcl1\gnl
  \gcl1\gev\gnl
  \gob1W\gend
  =
  \gbeg49
  \gvac3\got1V\gnl
  \gvac1\gdb\gcl1\gnl
  \gvac1\gcl1\gdnot{\pi_k}\glmptb\grmpt\gnl
  \gvac1\gsym\gnl
  \gdnot{\iota_k}\glmpb\grmptb\gcl2\gnl
  \gcl3\gcl2\gnl
  \gvac2\gdnot{\redu{(f_k)}}\glmptb\grmp\gnl
  \gcl1\gev\gnl
  \gob1W\gend
  =
  \gbeg49
  \gvac3\got1V\gnl
  \gvac1\gdb\gcl1\gnl
  \gvac1\gcl1\gdnot{\pi_k}\glmptb\grmpt\gnl
  \gvac1\gsym\gnl
  \gdnot{\iota_k}\glmpb\grmptb\gcl3\gnl
  \gcl3\gcl1\gnl
  \gvac1\gbmp{f_k}\gnl
  \gcl1\gev\gnl
  \gob1W\gend
  =
  \gbeg49
  \gvac3\got1V\gnl
  \gvac1\gdb\gcl1\gnl
  \gvac1\gcl1\gdnot{\pi_k}\glmptb\grmpt\gnl
  \gvac1\gsym\gnl
  \gdnot{\iota_k}\glmpb\grmptb\gcl3\gnl
  \gcl1\gcl1\gnl
  \gdnot{f_{k+1}}\glmptb\grmptb\gnl
  \gcl1\gev\gnl
  \gob1W\gend
  \]
  and thus, inductively
  \[f\FS_V^{(n)}
  =
  \gbeg49
  \gvac3\got1V\gnl
  \gvac1\gdb\gcl1\gnl
  \gvac1\gcl1\gdnot{\pi}\glmptb\grmpt\gnl
  \gvac1\gcn1113\gcn111{-1}\gnl
  \gdnot\iota\glmpb\grmptb\gcl1\gnl
  \gcl1\gev\gnl
  \gbmp f\gnl
  \gcl1\gnl
  \gob1W\gend
  =
  \gbeg 49
  \gvac3\got1V\gnl
  \gvac1\gdb\gcl1\gnl
  \gvac1\gcl1\gdnot{\pi_{0}}\glmptb\grmpt\gnl
  \gvac1\gsym\gnl
  \gdnot{\iota_{0}}\glmpb\grmptb\gcl3\gnl
  \gcl1\gcl1\gnl
  \gdnot{f_1}\glmptb\grmptb\gnl
  \gcl1\gev\gnl
  \gob1W\gend
  =
  \gbeg49
  \gvac3\got1V\gnl
  \gvac1\gdb\gcl1\gnl
  \gvac1\gcl1\gdnot{\pi_n}\glmptb\grmpt\gnl
  \gvac1\gsym\gnl
  \gdnot{\iota_n}\glmpb\grmptb\gcl3\gnl
  \gcl1\gcl1\gnl
  \gdnot{f_{n}}\glmptb\grmptb\gnl
  \gcl1\gev\gnl
  \gob1W\gend
  =
  \gbeg49
  \gvac3\got1V\gnl
  \gvac3\gcl1\gnl
  \gvac3\gbmp f\gnl
  \gvac1\gdb\gcl1\gnl
  \gvac1\gcl1\gdnot{\pi_n}\glmptb\grmpt\gnl
  \gvac1\gcn1113\gcn111{-1}\gnl
  \gdnot{\iota_n}\glmpb\grmptb\gcl1\gnl
  \gcl1\gev\gnl
  \gob1W\gend
  =
  \FS_W^{(n)}f {\quad.} \]
\end{proof}
\begin{cor}
  The Frobenius-Schur indicators are additive, that is
  $\nu_n(V\oplus W)=\nu_n(V)+\nu_n(W)$ holds for
  objects $V,W$ in a semisimple $k$-linear pivotal monoidal category $\CC$.
\end{cor}
\begin{proof}
Since the pivotal traces are cyclic, taking the pivotal trace is
additive (cf.\ \cite[Prop.1]{GelKaz:ITDM}), as is the natural
transformation $\FS^{(n)}$.
\end{proof}

\begin{remark}
  To define the Frobenius-Schur indicators and prove their basic
  properties, we need not assume that the category $\CC$ is
  semisimple. It is enough to assume that we are given a linear endofunctor
  $(\leer)^\triv$ of $\CC$ such
  that every $V^\triv$ is a direct sum of copies of the neutral
  object, and natural transformations
  $\iota\colon V^\triv\to V$ and $\pi\colon V\to V^\triv$.
  If $\CC=\C H$ is the category of finite-dimensional modules
  over a unimodular (quasi-)Hopf algebra $H$, we could take $V^\triv$ to be the
  submodule of $H$-invariants, and define $\pi$
  as multiplication by an integral.
  If $\iota$ and $\pi$ can be chosen so that
  $\pi\iota=\id$, then the triple $((\leer)^\triv,\iota,\pi)$ is
  unique up to isomorphism and the resulting Frobenius-Schur
  endomorphisms are unique. In the case that $\CC=\C H$ as above, the existence of
  such a triple is well-known to imply semisimplicity of $\CC$.
\end{remark}

It is also possible to define Frobenius-Schur endomorphisms whose
traces correspond to the indicators $\nu_{n,k}(V)$ with $k>1$.
Since it is quite tedious to even write these down in full
generality in the non-strict case, we will give details only in
the strict case.
\begin{defn}Let $V$ be an object in a left and right rigid monoidal
category.
\begin{enumerate}
\item If $\CC$ is strict monoidal, we define generalized
Frobenius-Schur endomorphisms of $V$ by
\[\FS_V^{(n,\ell,r)}=
  \gbeg {10}9
  \gvac7\got1V\gnl
  \gwdb7\gcl2\gnl
  \gcl4\gvac3\gdb\gcl1\gvac1\gdb\gnl
  \gvac4\gcl1\glmptb\gdnot\pi\gcmpt\gcmpt\grmpt\gcl4\gnl
  \gvac4\gsym\gnl
  \gvac1\glmpb\gdnot\iota\gcmpb\gcmpb\grmptb\gcl1\gnl
  \gev\gcl2\gcl1\gev\gnl
  \gvac3\gwev7\gnl
  \gvac2\gob1V\gend\]
where we have used the following distribution of the tensor
factors over the ``legs'' of the graphical symbols for $\pi$ and
$\iota$, with $k=\ell+r+1$:
\[\gbeg{10}4
  \got3{V^{\o(n-k)}}\got3{V^{\o\ell}}\got1V\got3{V^{\o r}}\gnl
  \gvac1\gcl1\gvac2\gcl1\gvac1\gcl1\gvac1\gcl1\gnl
  \gvac1\glmpt\gcmp\gcmp\gdnot\pi\gcmpt\gcmp\gcmpt\gcmp\grmpt\gnl
  \gend
  \qquad
  \text{ and }
  \qquad
  \gbeg{10}4
  \gnl
  \gvac1\glmpb\gcmp\gcmpb\gdnot{\iota}\gcmp\gcmpb\gcmp\gcmp\grmpb\gnl
  \gvac1\gcl1\gvac1\gcl1\gvac1\gcl1\gvac2\gcl1\gnl
  \gob3{V^{\o\ell}}\gob1V\gob3{V^{\o r}}\gob3{V^{\o{(n-k)}}}\gend\]
  We abbreviate $\FS_V^{(n,k)}=\FS_V^{(n,k-1,0)}$ for $1\leq k<n$.
\item In general, define
  $\FS_V^{(n,\ell,r)}=\FF\inv\left(\FS_{\FF(V)}^{(n,\ell,r)}\right)$, where
  $\FF\colon\CC\to\DD$ is a monoidal equivalence with a strict
  monoidal category $\DD$.
\end{enumerate}
\end{defn}
\begin{remark}
  Strictly speaking, the second part of the definition requires
  the proof that the endomorphisms defined in the first part are
  invariant under monoidal equivalences between strict monoidal
  categories. Except for the difficulty to write down the large
  expressions, this is routine, however. Of course the
  general definition is then also invariant under monoidal category
  equivalences. On these grounds, we will only give proofs for the
  strict case below.
\end{remark}
\begin{prop}
Let $V$ be an object of a $k$-linear semisimple pivotal monoidal
category.

We have $\nu_{n,k}(V)=\ptrl\left(\FS_V^{(n,k)}\right)$ for all
$n,k$.

If $r>0$, we have
$\ptrl\left(\FS_V^{(n,\ell,r)}\right)=\ptrr\left(\FS_V^{(n,\ell+1,r-1)}\right)$.

If $\CC$ is spherical, we have
$\nu_{n,k}(V)=\ptr\left(\FS_V^{(n,\ell,r)}\right)$ whenever
$k=\ell+r+1$.
\end{prop}
\begin{proof}
The proof of the first assertion is analogous to that of
\nmref{traceformula}. The second is obvious from the graphical
representation of left and right pivotal traces, and the third is
a consequence of the first two.
\end{proof}
\begin{remark}\nmlabel{Remark}{FSsame}
If $V$ is a simple object in $\CC$, then $\ptrl(\id_V)\neq
0\neq\ptrr(\id_V)$. Indeed $\ptrl(\id_V)$, for example, is the
composition of two nonzero morphisms $I\to V\o V\du$ and $V\o
V\du\to I$, which is nonzero since the multiplicity of $I$ in $V\o
V\du$ is one; see \cite[Lem.1]{GelKaz:ITDM} or the discussion
after \cite[Prop.2.1]{EtiNikOst:FC}.

In particular, an endomorphism of $V$ is uniquely determined by
its trace. Thus the preceding proposition shows that
$\FS_V^{(n,\ell,r)}=\FS_V^{(n,\ell+r+1)}$ holds for every simple
object $V$ in a $k$-linear semisimple spherical monoidal category.
Since the next proposition shows that the various Frobenius-Schur
endomorphisms are natural,
$\FS_V^{(n,\ell,r)}=\FS_V^{(n,\ell+r+1)}$ then holds for all
objects $V$.
\end{remark}

\begin{prop}
If $n$ and $k=\ell+r+1$ are relatively prime, then
$\FS_V^{(n,\ell,r)}$ is natural in $V$.
\end{prop}
\begin{proof}
We use the permutation $s\in S_n$ determined by requiring
$s(i)\in\{1,\dots,n\}$ to be congruent to $i+k$ modulo $n$.
(This occurs for similar reasons in  \cite{KasSomZhu:HFSI} as a
special case of more general permutation constructions.) Note that
$s$ is an $n$-cycle since $n$ and $k$ are relatively prime.

Consider $f\colon V\to W$. For any $X_i,Y_j$ in $\CC$ we will
write
 \[f_p
 \colon X_1\o\dots\o X_{p-1}\o V\o Y_1\o\dots\o Y_u
 \to X_1\o\dots\o X_{p-1}\o W\o Y_1\o\dots\o Y_u.\]
for the morphism that acts as $f$ in the $p$-th position and the
identity otherwise. Define a series of objects and morphisms
\[V[0]\xrightarrow{f[1]}V[1]\xrightarrow{f[2]}V[2]\to\dots\xrightarrow{f[n]}V[n]\]
by $V[0]=V^{\o n}$ and $f[i]=f_{s^{i-1}(\ell+1)}$; this fixes
$V[i]$, which has to be the appropriate target; note that the
sequence is well-defined since $s$ is transitive; we have
$V[n]=W^{\o n}$. Denote the injection from and projection to the
trivial component of $V[i]$ by $\iota_i$ and $\pi_i$.

Now as before we have
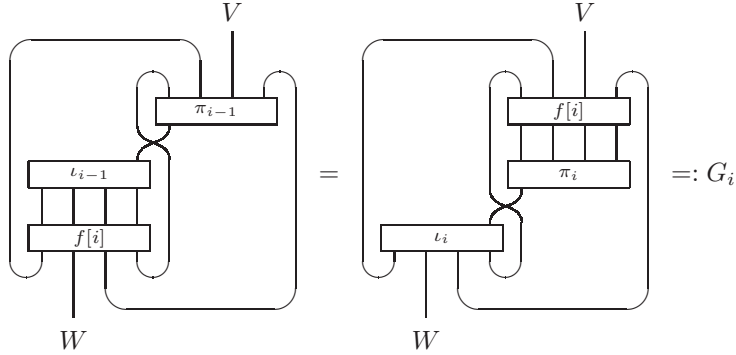
\[\gbeg {10}{11}
  \gvac7\got1V\gnl
  \gwdb7\gcl2\gnl
  \gcl6\gvac3\gdb\gcl1\gvac1\gdb\gnl
  \gvac4\gcl1\glmptb\gdnot{\pi_{i-1}}\gcmpt\gcmpt\grmpt\gcl6\gnl
  \gvac4\gsym\gnl
  \gvac1\glmpb\gdnot{\iota_{i-1}}\gcmpb\gcmpb\grmptb\gcl3\gnl
  \gvac1\gcl1\gcl1\gcl1\gcl1\gnl
  \gvac1\glmptb\gdnot{f[i]}\gcmptb\gcmptb\grmptb\gnl
  \gev\gcl2\gcl1\gev\gnl
  \gvac3\gwev7\gnl
  \gvac2\gob1W\gend
  =
  \gbeg {10}{11}
  \gvac7\got1V\gnl
  \gwdb7\gcl2\gnl
  \gcl6\gvac3\gdb\gcl1\gvac1\gdb\gnl
  \gvac4\gcl3\glmptb\gdnot{f[i]}\gcmptb\gcmptb\grmptb\gcl6\gnl
  \gvac5\gcl1\gcl1\gcl1\gcl1\gnl
  \gvac4\gvac1\glmptb\gdnot{\pi_i}\gcmpt\gcmpt\grmpt\gnl
  \gvac4\gsym\gnl
  \gvac1\glmpb\gdnot{\iota_i}\gcmpb\gcmpb\grmptb\gcl1\gnl
  \gev\gcl2\gcl1\gev\gnl
  \gvac3\gwev7\gnl
  \gvac2\gob1W\gend
  =:G_i\]
  for any $i\in\{1,\dots,n\}.$
  If $p:=s^{i-1}(\ell+1)\leq n-k$, then
  \[G_i
  =
  \gbeg {10}{11}
  \gvac7\got1V\gnl
  \gwdb7\gcl2\gnl
  \gcl6\gvac3\gdb\gcl1\gvac1\gdb\gnl
  \gvac4\gcl3\gbmp{f_{p}}\gcl1\gcl1\gcl1\gcl6\gnl
  \gvac5\gcl1\gcl1\gcl1\gcl1\gnl
  \gvac4\gvac1\glmptb\gdnot{\pi_i}\gcmpt\gcmpt\grmpt\gnl
  \gvac4\gsym\gnl
  \gvac1\glmpb\gdnot{\iota_i}\gcmpb\gcmpb\grmptb\gcl1\gnl
  \gev\gcl2\gcl1\gev\gnl
  \gvac3\gwev7\gnl
  \gvac2\gob1W\gend
  =
  \gbeg {10}{11}
  \gvac7\got1V\gnl
  \gwdb7\gcl2\gnl
  \gcl6\gvac3\gdb\gcl1\gvac1\gdb\gnl
  \gvac4\gcl1\glmptb\gdnot{\pi_{i}}\gcmpt\gcmpt\grmpt\gcl6\gnl
  \gvac4\gsym\gnl
  \gvac1\glmpb\gdnot{\iota_{i}}\gcmpb\gcmpb\grmptb\gcl3\gnl
  \gvac1\gcl1\gcl1\gcl1\gcl1\gnl
  \gvac1\gcl1\gcl1\gcl1\gbmp{f_{p}}\gnl
  \gev\gcl2\gcl1\gev\gnl
  \gvac3\gwev7\gnl
  \gvac2\gob1W\gend
  =
  \gbeg {10}{11}
  \gvac7\got1V\gnl
  \gwdb7\gcl2\gnl
  \gcl6\gvac3\gdb\gcl1\gvac1\gdb\gnl
  \gvac4\gcl1\glmptb\gdnot{\pi_{i}}\gcmpt\gcmpt\grmpt\gcl6\gnl
  \gvac4\gsym\gnl
  \gvac1\glmpb\gdnot{\iota_{i}}\gcmpb\gcmpb\grmptb\gcl3\gnl
  \gvac1\gcl1\gcl1\gcl1\gcl1\gnl
  \gvac1\glmptb\gdnot{f[i+1]}\gcmptb\gcmptb\grmptb\gnl
  \gev\gcl2\gcl1\gev\gnl
  \gvac3\gwev7\gnl
  \gvac2\gob1W\gend\]
  since $p+k=s(p)=s^{i}(\ell+1)$.
  If $n-k<p\leq n-k+\ell$, then for $q:=p-n+k$
  \[G_i
  =
  \gbeg {10}{11}
  \gvac7\got1V\gnl
  \gwdb7\gcl2\gnl
  \gcl6\gvac3\gdb\gcl1\gvac1\gdb\gnl
  \gvac4\gcl3\gcl1\gbmp{f_{q}}\gcl1\gcl1\gcl6\gnl
  \gvac5\gcl1\gcl1\gcl1\gcl1\gnl
  \gvac4\gvac1\glmptb\gdnot{\pi_i}\gcmpt\gcmpt\grmpt\gnl
  \gvac4\gsym\gnl
  \gvac1\glmpb\gdnot{\iota_i}\gcmpb\gcmpb\grmptb\gcl1\gnl
  \gev\gcl2\gcl1\gev\gnl
  \gvac3\gwev7\gnl
  \gvac2\gob1W\gend
  =
  \gbeg {10}{11}
  \gvac7\got1V\gnl
  \gwdb7\gcl2\gnl
  \gcl6\gvac3\gdb\gcl1\gvac1\gdb\gnl
  \gvac4\gcl1\glmptb\gdnot{\pi_{i}}\gcmpt\gcmpt\grmpt\gcl6\gnl
  \gvac4\gsym\gnl
  \gvac1\glmpb\gdnot{\iota_{i}}\gcmpb\gcmpb\grmptb\gcl3\gnl
  \gvac1\gcl1\gcl1\gcl1\gcl1\gnl
  \gvac1\gbmp{f_{q}}\gcl1\gcl1\gcl1\gnl
  \gev\gcl2\gcl1\gev\gnl
  \gvac3\gwev7\gnl
  \gvac2\gob1W\gend
  =
  \gbeg {10}{11}
  \gvac7\got1V\gnl
  \gwdb7\gcl2\gnl
  \gcl6\gvac3\gdb\gcl1\gvac1\gdb\gnl
  \gvac4\gcl1\glmptb\gdnot{\pi_{i}}\gcmpt\gcmpt\grmpt\gcl6\gnl
  \gvac4\gsym\gnl
  \gvac1\glmpb\gdnot{\iota_{i}}\gcmpb\gcmpb\grmptb\gcl3\gnl
  \gvac1\gcl1\gcl1\gcl1\gcl1\gnl
  \gvac1\glmptb\gdnot{f[i+1]}\gcmptb\gcmptb\grmptb\gnl
  \gev\gcl2\gcl1\gev\gnl
  \gvac3\gwev7\gnl
  \gvac2\gob1W\gend\]
  since $q=s(p)=s^{i}(\ell+1)$.  If $p>n-k+\ell+1=n-r+1$, then we set
  $q=p-n+k-\ell-1$ and find
  \[G_i=
  \gbeg {10}{11}
  \gvac7\got1V\gnl
  \gwdb7\gcl2\gnl
  \gcl6\gvac3\gdb\gcl1\gvac1\gdb\gnl
  \gvac4\gcl3\gcl1\gcl1\gcl1\gbmp{f_{q}}\gcl6\gnl
  \gvac5\gcl1\gcl1\gcl1\gcl1\gnl
  \gvac4\gvac1\glmptb\gdnot{\pi_i}\gcmpt\gcmpt\grmpt\gnl
  \gvac4\gsym\gnl
  \gvac1\glmpb\gdnot{\iota_i}\gcmpb\gcmpb\grmptb\gcl1\gnl
  \gev\gcl2\gcl1\gev\gnl
  \gvac3\gwev7\gnl
  \gvac2\gob1W\gend
  =
  \gbeg {10}{11}
  \gvac7\got1V\gnl
  \gwdb7\gcl2\gnl
  \gcl6\gvac3\gdb\gcl1\gvac1\gdb\gnl
  \gvac4\gcl1\glmptb\gdnot{\pi_{i}}\gcmpt\gcmpt\grmpt\gcl6\gnl
  \gvac4\gsym\gnl
  \gvac1\glmpb\gdnot{\iota_{i}}\gcmpb\gcmpb\grmptb\gcl3\gnl
  \gvac1\gcl1\gcl1\gcl1\gcl1\gnl
  \gvac1\gcl1\gcl1\gbmp{f_{q}}\gcl1\gnl
  \gev\gcl2\gcl1\gev\gnl
  \gvac3\gwev7\gnl
  \gvac2\gob1W\gend
  =
  \gbeg {10}{11}
  \gvac7\got1V\gnl
  \gwdb7\gcl2\gnl
  \gcl6\gvac3\gdb\gcl1\gvac1\gdb\gnl
  \gvac4\gcl1\glmptb\gdnot{\pi_{i}}\gcmpt\gcmpt\grmpt\gcl6\gnl
  \gvac4\gsym\gnl
  \gvac1\glmpb\gdnot{\iota_{i}}\gcmpb\gcmpb\grmptb\gcl3\gnl
  \gvac1\gcl1\gcl1\gcl1\gcl1\gnl
  \gvac1\glmptb\gdnot{f[i+1]}\gcmptb\gcmptb\grmptb\gnl
  \gev\gcl2\gcl1\gev\gnl
  \gvac3\gwev7\gnl
  \gvac2\gob1W\gend\]
  since now $q+\ell+1=p+k-n=s(p)=s^{i}(\ell+1)$. Because $s$ has
  order $n$ and $s(n-k+\ell+1)=\ell+1$, the case $p=n-k+\ell+1=s^{n-1}(\ell+1)$
  occurs when $i=n$. Thus, using induction, we have
  \[f\FS^{(n,\ell,r)}_V
  =\gbeg {10}{11}
  \gvac7\got1V\gnl
  \gwdb7\gcl2\gnl
  \gcl6\gvac3\gdb\gcl1\gvac1\gdb\gnl
  \gvac4\gcl1\glmptb\gdnot{\pi_{0}}\gcmpt\gcmpt\grmpt\gcl6\gnl
  \gvac4\gsym\gnl
  \gvac1\glmpb\gdnot{\iota_{0}}\gcmpb\gcmpb\grmptb\gcl3\gnl
  \gvac1\gcl1\gcl1\gcl1\gcl1\gnl
  \gvac1\glmptb\gdnot{f[1]}\gcmptb\gcmptb\grmptb\gnl
  \gev\gcl2\gcl1\gev\gnl
  \gvac3\gwev7\gnl
  \gvac2\gob1W\gend
  =
  \gbeg {10}{11}
  \gvac7\got1V\gnl
  \gwdb7\gcl2\gnl
  \gcl6\gvac3\gdb\gcl1\gvac1\gdb\gnl
  \gvac4\gcl3\glmptb\gdnot{f[n]}\gcmptb\gcmptb\grmptb\gcl6\gnl
  \gvac5\gcl1\gcl1\gcl1\gcl1\gnl
  \gvac4\gvac1\glmptb\gdnot{\pi_n}\gcmpt\gcmpt\grmpt\gnl
  \gvac4\gsym\gnl
  \gvac1\glmpb\gdnot{\iota_n}\gcmpb\gcmpb\grmptb\gcl1\gnl
  \gev\gcl2\gcl1\gev\gnl
  \gvac3\gwev7\gnl
  \gvac2\gob1W\gend
  =
  \FS^{(n,\ell,r)}_Wf
  \]
\end{proof}

\begin{remark}
As our definition of $E_{V}^{(n)}$, the definition of the various
Frobenius-Schur endomorphisms above is not left-right symmetric.
In particular, the Frobenius-Schur endomorphism $\FS^{(n)}$ for
$V^\sym$ as an object of $\CC^\sym$ is given pictorially by
\[\gbeg47
  \got1V\gnl
  \gcl1\gdb\gnl
  \gdnot\pi\glmpt\grmptb\gcl1\gnl
  \gvac1\gsym\gnl
  \gvac1\gcl1\gdnot\iota\glmptb\grmpb\gnl
  \gvac1\gev\gcl1\gnl
  \gvac3\gob1V\gend\]
and similarly for $\FS^{(n,k)}$. Using \nmref{symmetric} we see that
(some of) these endomorphisms could equally well be used to compute
the indicators. Moreover, as in \nmref{FSsame}, for the case of a
$k$-linear semisimple spherical monoidal category the endomorphisms
agree up to relabelling, e.g.\ the one instance for which we have
given a picture agrees with $\FS_V^{(n,n-1)}$.
\end{remark}
\bibliographystyle{acm}
%\bibliography{eigene,andere}
%\end{document}
\def\cprime{$'$}

\end{document}